\documentclass{article}
\usepackage[utf8]{inputenc}
\usepackage{lineno,epsfig,graphics,amssymb,amsmath,graphicx,pifont,soul,esint,float,multirow,array,color,url,xcolor,booktabs,mathrsfs,mathdots}
\usepackage{siunitx}
\usepackage{enumitem}
\usepackage[margin=1in]{geometry}
\usepackage{url}
\usepackage{caption}
\usepackage{subcaption}
\usepackage{cleveref}
\usepackage{mathtools}
\usepackage{authblk}

\newcommand{\bl}[1]{\boldsymbol{#1}}
\newcommand{\cola}[1]{\textcolor{black}{#1}}

\title{Stabilized finite element methods for the time-spectral convection-diffusion equation} 
\author[1,*]{Mahdi Esmaily}
\author[1]{Dongjie Jia}
\affil[1]{Sibley School of Mechanical and Aerospace Engineering, Cornell University, Ithaca, NY, USA}
\affil[*]{Correspondence: me399@cornell.edu}

\begin{document}
\maketitle

\begin{abstract}
     Discretizing a solution in the Fourier domain rather than the time domain presents a significant advantage in solving transport problems that vary smoothly and periodically in time, such as cardiorespiratory flows.
     The finite element solution of the resulting time-spectral formulation is investigated here for the convection-diffusion equations. 
     In addition to the baseline Galerkin's method, we consider stabilized approaches inspired by the streamline upwind Petrov/Galerkin (SUPG), Galerkin/least square (GLS), and variational multiscale (VMS) methods. 
     We also introduce a new augmented SUPG (ASU) method that, by design, produces a nodally exact solution in one dimension for piecewise linear interpolation functions. 
     Comparing these five methods using 1D, 2D, and 3D canonical test cases shows while the ASU is most accurate overall, it exhibits stability issues in extremely oscillatory flows with a high Womersley number in 3D. 
     The GLS method, which is identical to the VMS for this problem, presents an attractive alternative due to its excellent stability and reasonable accuracy. 
\end{abstract}

\section{Introduction}
In recent decades, the use of cardiorespiratory simulations for surgical design~\cite{Esmaily2015ABG, Marsden20081890}, diagnosis~\cite{driessen2019comparison}, and patient-specific modeling~\cite{taylor2009patient} has been on the rise. 
The finite element method has emerged as a popular choice among various numerical methods due to its ability to handle complex physics and geometries effectively~\cite{taylor1998finite, Vignone20103, Bazilevs2009computational}. 
However, the computational cost associated with these simulations, particularly for inverse problems like optimization~\cite{Yang20102135, Esmaily2012optimization, verma2018optimization}, parameter identification~\cite{tran2017automated, mirramezani2020distributed, pfaller2021automated}, and uncertainty quantification~\cite{phillips2003quantifying, sankaran2011stochastic, schiavazzi2017generalized} has limited expansion of their use beyond academic settings and to, for instance, industrial applications.

The cost of a computational fluid dynamics simulation is determined by the dimension of the discrete problem, which is the product of the degrees of freedom for spatial discretization (i.e., number of grid points) and the number of time steps for time integration. 
In time-dependent simulations, the number of time steps can be quite large, potentially reaching thousands as the time step size should be small enough to ensure accuracy and the time integration period should be long enough to ensure the independence of the solution from the arbitrary initial condition~\cite{arbia2014numerical}. 
Given this large number, the cost of a simulation can be dramatically reduced if one were to replace the time integration with a more cost-effective alternative.

For this purpose, we propose discretizing the fluid problem in the frequency rather than time domain~\cite{hupp2016parallel, hall2002computation, arbenz2017comparison}. 
This choice is motivated by the fact that the transport variables in the cardiorespiratory system often vary smoothly and periodically in time. 
As a result, those variables can be well-approximated with only a handful of Fourier modes~\cite{meng2020complex,esmaily2022stabilized}. 
Having a handful of modes, in contrast to thousands of time steps, presents a huge opportunity for reducing the dimensionality of the discrete problem and as a consequence, the cost of a typical cardiorespiratory simulation. 

To realize this cost advantage, a robust and efficient time-spectral incompressible Navier-Stokes equations solver must be developed. 
Doing so poses three main challenges: dealing with the incompressibility constraint, ensuring stability under strongly convective regimes, and efficiently handling mode coupling associated with the nonlinear convective acceleration term, which we leave for future studies. 

The first challenge above, which can be replicated exclusively in the unsteady Stokes equation, has been the subject of our earlier studies~\cite{meng2020time,esmaily2022stabilized}. 
In the first study~\cite{meng2020time}, we showcased the possibility of reducing the cost of simulation by several orders of magnitude if one were to solve it in the frequency domain.
The later study~\cite{esmaily2022stabilized} generalized that method to avoid complex arithmetic and circumvent the inf-sup condition to allow for the use of similar interpolation functions for pressure and velocity~\cite{ladyzhenskaya1969mathematical, babuvska1971error, brezzi1974existence}. 

The focus of the present study is the second challenge: dealing with strongly convective flows. Specifically, we aim to identify a finite element method that produces accurate and stable solutions for a broad range of flow conditions. To decouple this challenge from the third mode-coupling challenge, we concentrate on the linear unsteady convection-diffusion equation driven by a given steady flow. Additionally, we limit our discussion to linear interpolation functions due to their widespread use in the discretization of complex cardiorespiratory geometries.

In the context of conventional time formulation, it is well-known that Galerkin's method produces nonphysical oscillations in strongly convective regimes~\cite{johnson2012numerical}. 
The literature dedicated to dealing with this issue is too vast to recount here~\cite{hughes1987recent, franca1992stabilized, bazilevs2013computational, hughes2018multiscale}. 
Thus, we forgo discussion of methods such as discontinuity capturing~\cite{Hughes1986beyond, hughes1986discontinuity, tezduyar1986discontinuity, takizawa2018stabilization} and bubble functions \cite{brezzi1998further, franca2000improved, franca2002stability} and focus on a select few popular stabilization methods that are investigated in details in the following sections. 

One of the earliest techniques proposed to counteract instabilities associated with convection-dominated flows was to use an upwinding scheme by appropriately adjusting the test function weights in the upstream and downstream of the tested node~\cite{hughes1979multidimentional, brooks1982streamline}.
This strategy, known as the streamline upwind Petrov/Galerkin (SUPG) method, successfully generates stable and accurate results for a conventional time formulation by adding a direction-dependent diffusion to the underlying Galerkin's formulation.
This study investigates the extent to which the issues with the Galerkin's method carry over to the time-spectral form of the convection-diffusion equation and evaluates the effectiveness of SUPG in eliminating those issues.

Additionally, we explore other stabilization methods suitable for convection-dominated flows, including the Galerkin/least-squares-based method (GLS) that introduces a symmetric residual penalty term to the discrete form~\cite{hughes1989new, shakib1991new, shakib1989finite} and the variational multiscale (VMS) method that models the unresolved scales in the discrete solution through its residual~\cite{hughes1995multiscale, hughes2007variational, bazilevs2007variational}.

In addition to the above stabilized methods, which are adaptations of conventional time formulations, we propose a new method designed specifically for the time-spectral form of the convection-diffusion equation. This new method, which can be viewed as an augmented SUPG method (ASU), is designed to accomplish what the SUPG is designed to accomplish, but for the time-spectral rather than the steady-state version of the convection-diffusion equation. More specifically, we design the ASU to produce a nodally exact solution for a one-dimensional model problem for the time-spectral convection-diffusion equation, regardless of whether the solution is steady or unsteady.

Beside convection-dominated flows, the challenges associated with the solution of unsteady Stokes equations in the frequency domain share common ground with another body of literature that concerns the time-harmonic solution of Helmholtz equation to model acoustics \cite{harari2006survey, harari1990design, harari1992cost}. 
A significant hurdle in computational acoustics revolves around accurately representing waves with short wavelengths relative to the mesh size. If left unattended, this issue can pollute the solution and introduce dispersive errors \cite{babuvska1995generalized, babuska1997pollution}. 
The presence of this problem, which also surfaces in scenarios characterized by high Womersley numbers in the present case, has garnered considerable attention in various studies \cite{franca2000improved, franca1997residual}, leading to methods inspired by the VMS and GLS to alleviate the above deficiencies \cite{thompson1995galerkin, oberai1998multiscale, oberai2000residual, harari1992galerkin}. 
By evaluating the performance of these methods in regimes of weak convection but strong oscillation, the findings of this study hold relevance for individuals interested in employing finite element techniques for time-harmonic acoustics.

The article is structured as follows: In Section \ref{sec:1d}, we introduce a one-dimensional model problem to derive the ASU method and discuss the extension of other methods to time-spectral form. Section \ref{sec:md} discusses the generalization of these methods to multiple dimensions, and Section \ref{sec:conclusions} concludes the study.

\section{A 1D model problem}\label{sec:1d}
To rigorously develop an accurate technique for solving the convection-diffusion problem in multiple dimensions, we first turn to a simple one-dimensional problem. 
This 1D problem has historical significance for its roots in the development of the SUPG method. 
We also use this model problem to construct the ASU approach for the time-spectral convection-diffusion equations. 

\subsection{Problem statement} \label{sec:1D-problem-statement}
Consider the unsteady convection of a neutral tracer $\hat \phi(x,t)$ in a one-dimensional domain that is governed by
\begin{equation}
\begin{split}
    \hat \phi_{,t} + a\hat \phi_{,x} &= \kappa \hat \phi_{,xx}, \\
    \hat \phi(0,t) &= 0, \\
    \hat \phi(L,t) &= \cos(\omega t),
\end{split}
\label{1D-time}
\end{equation}
where $L$ is the domain size, $\kappa \in \mathbb R^+$ is the diffusivity, $\omega\in \mathbb R$ is the oscillation frequency of the boundary condition, and $a\in \mathbb R$ is the convective velocity that is uniform in the entire domain.
The initial condition was not specified in Eq. \eqref{1D-time} because we are solely interested in its particular solution (i.e., $\hat \phi(x,t)$ as $t \to \infty$) that is independent of the initial transient behavior of $\hat \phi$ when $\kappa > 0$. 

Even though the boundary condition specified in Eq. \eqref{1D-time} is expressed in the form of a uni-modal excitation, one can simply generalize what we discuss below to any arbitrary (but well-behaved) time-varying boundary condition $f(t)$ given that $\hat \phi$ is linear in terms of $f$ and $f$ can be expressed as the summation of trigonometric functions through Fourier transformation.

Since we are interested in the time-spectral formulation of Eq. \eqref{1D-time}, we instead attempt to solve an equivalent problem that is
\begin{equation}
\begin{split}
    i\omega \phi + a\phi_{,x} &= \kappa \phi_{,xx}, \\
    \phi(0) &= 0, \\
    \phi(L) &= 1,
\end{split}
\label{1D-fourier}
\end{equation}
where $i^2=-1$ and
\begin{equation}
    \hat \phi(x,t) = {\rm Real}\left(\phi(x)e^{i\omega t}\right) = \phi_r \cos(\omega t) - \phi_i \sin(\omega t). 
    \label{1D-FT}
\end{equation} 
In Eq. \eqref{1D-FT}, $\phi_r = \phi_r(x)$ and $\phi_i = \phi_i(x)$ denote the real and imaginary components of $\phi$. 
These two functions determine the overall amplitude $|\phi| = \sqrt{\phi_r^2 + \phi_i^2}$ of the solution and its phase shift $\theta = \tan^{-1}(\phi_i/\phi_r)$ relative to the boundary condition. 
Since $\phi_r$ and $\phi_i$ capture the overall behavior of the solution $\hat \phi(x,t)$, we rely on these two functions to evaluate the various methods below.  

\subsection{Exact solution}
Since Eq. \eqref{1D-fourier} is a constant coefficient second-order ordinary differential equation, its solution can be obtained through elementary means and is
\begin{equation}
    \phi(x) =\frac{\exp\left(r_1 \frac{x}{L}\right)-\exp\left(r_2 \frac{x}{L}\right)}{\exp(r_1)-\exp(r_2)},
    \label{1D-exact}
\end{equation}
where
\begin{equation}
    r_{1,2} = P \pm \sqrt{P^2+ iW^2},
\label{r12}
\end{equation}
are the roots of the characteristic polynomial and are expressed in terms of 
\begin{align}
P &= \frac{aL}{2\kappa}, \label{P-def}\\
W &= L\sqrt{\frac{\omega}{\kappa}},\label{W-def}
\end{align}
which are the Peclet and Womersley numbers, respectively. 

The Peclet number, which weighs convection relative to the diffusion, is analogous to the Reynolds number and is typically much larger than one for engineering applications or blood flow in large blood vessels. 
The Womersley number, on the other hand, measures the importance of unsteady effects against diffusion and can vary from tens in major vessels to values much smaller than one in the cardiovascular system~\cite{Womersley1955method,feintuch2007hemodynamics}. 

\subsection{Baseline Galerkin's method}
Galerkin's approximate solution is obtained from the weak form of Eq. \eqref{1D-fourier} by multiplying it by a test function $w^h(x)$ and integrating by parts the diffusion while noting the approximate $\phi^h$ diminishes at the boundaries. 
The resulting problem statement is to find $\phi^h(x)$ that satisfied the boundary conditions such that for any $w^h(x)$, which diminishes at the boundaries, we have
\begin{equation}
    (w^h,i\omega \phi^h) + (w^h,a\phi^h_{,x}) + (w^h_{,x},\kappa \phi^h_{,x}) = 0, 
\label{1D-weak}
\end{equation}
where $(f,g) = \int_0^Lfg {\rm d}x$ denotes the inner product of functions $f$ and $g$.  

As detailed in Appendix \ref{appendix1}, Eq. \eqref{1D-weak} has a closed-form solution for piecewise linear interpolation functions of uniform size $h$. 
Taking $N$ to be the number of elements (so that the nodal position $x_A=hA$ for $A=0,1,\cdots,N$), the solution at node $A$ is
\begin{equation}
    \phi^h(x_A) = \frac{\rho_1^A - \rho_2^A}{\rho_1^N - \rho_2^N},
    \label{1d-gal}
\end{equation}
where 
\begin{equation}
    \rho_{1,2} = \frac{1 + 2i\beta \pm \sqrt{\alpha^2 - 3\beta^2 + 6i\beta}}{1 - \alpha - i\beta},
\label{rho12} 
\end{equation}
with 
\begin{align}
\alpha &= \frac{ah}{2\kappa}, 
\label{alpha}\\
\beta &= \frac{\omega h^2}{6\kappa}. 
\label{beta}
\end{align}

The two variables $\alpha$ and $\beta$ in Eqs. \eqref{alpha} and \eqref{beta} are the Peclet number and square of the Womersley number, respectively, that are defined based on the element size $h$ rather than the domain size $L$. 
Differently put, $\alpha$ represents the relative magnitude of the convective term in comparison to the diffusive term at the element length scale. 
Similarly, $\beta$ represents the relative magnitude of the acceleration term in comparison to the diffusive term at the element length scale. 

In the steady state flows, where $\omega=0$ and thus $\beta=0$, $\rho_{1,2}=1,(1+\alpha)/(1-\alpha)$. 
This result, which has been established in the past \cite{brooks1982streamline}, explains the nonphysical oscillatory nature of Galerkin's solution in strongly convective flows in which $|\alpha|>1$. 
In these cases $\rho_2<0$ will produce alternating signs for $\phi^h(x_A)$ as $A$ switches between odd and even numbers. 
In the next section, we will show this behavior persists in unsteady flows when $\beta> 0$. 

\subsection{The issue} \label{sec:issue}
As stated earlier, it is well established that the spatiotemporal Galerkin's method fails at strongly convective regimes by generating nonphysical oscillations in the solution. 
To what extent this issue persists in unsteady regimes is what we investigate below. 

To answer the above question, we considered the 1D model problem stated in Eq. \eqref{1D-fourier}. 
We prescribe a flow that is from right to left ($a<0$) so that the solution temporal oscillations (that are physical and caused by the unsteady boundary conditions) propagate into the computational domain.
This problem setup allows us to better contrast various methods. 
Although results are not presented here for $a>0$, the error for those cases resembles what is presented below for $a<0$.  

\begin{figure}
  \centering
  \include{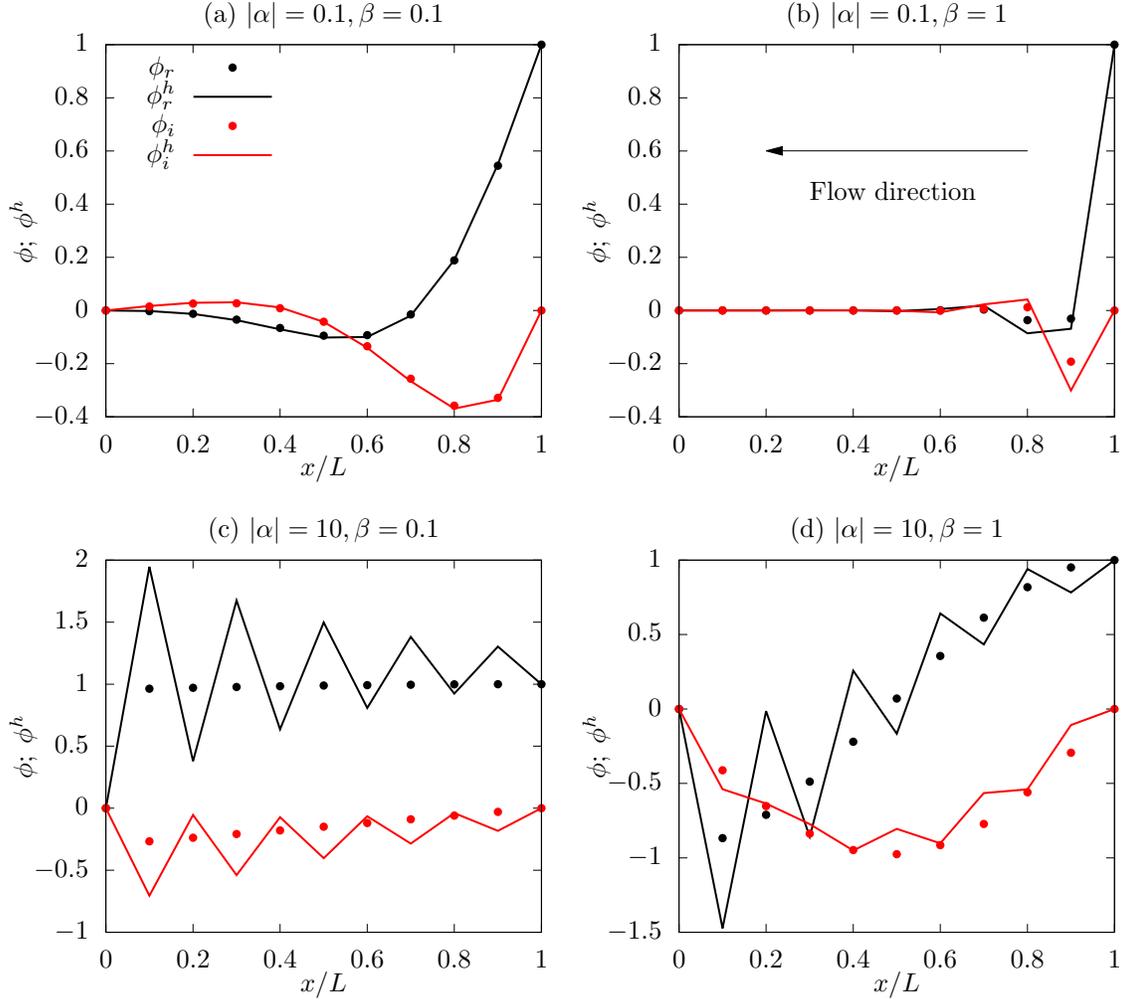}
   \caption{The real (black) and imaginary (red) part of the solution to the 1D time-spectral convection-diffusion problem. 
   The exact (dots) and Galerkin's (solid line) solutions are graphed for $|\alpha|=0.1$ (top row), $|\alpha|=10$ (bottom row), $\beta=0.1$ (left column), and $\beta=1$ (right column).}
  \label{fig:issue}    
\end{figure}

The real and imaginary parts of the solution obtained from Galerkin's method are compared against the exact solution from Eq. \eqref{1D-exact} in Figure \ref{fig:issue}  
The results are obtained for a wide range of element Peclet and Womersley numbers to demonstrate three key observations.
Firstly, in strongly convective flows (large $\alpha$), the error is dominated by nonphysical oscillations (Figure \ref{fig:issue}-(c,d)). 
These oscillations, which are similar to that of the steady problem, are generated due to the presence of sharp changes near the Dirichlet boundary. 
Secondly, at relatively small $\alpha$ and $\beta$, which correspond to cases where the mesh is sufficiently small, Galerkin's method provides a very good approximation of the solution (Figure \ref{fig:issue}-(a)). 
Lastly, in highly oscillatory but weakly convective flows Galerkin's solution overshoots the exact solution near the oscillatory boundary (Figure \ref{fig:issue}-(b)).
In the next section, we will discuss how various stabilization techniques overcome these issues. 

\subsection{Streamline upwind Petrov/Galerkin method}\label{sec:1dsupg}
The SUPG method is constructed by modifying the way that the steady state version of Eq. \eqref{1D-fourier} is tested~\cite{brooks1982streamline}. 
Instead of testing the equation with $w$ as shown in Eq. \eqref{1D-weak}, the SUPG adds an upwinding contribution and tests it with $w + \tau a w_{,x}$.  
The added term $\tau a w_{,x}$, which is only active along the streamwise direction in multidimensional flows, increases the overall test function weight upstream of the tested point. 
With this adjustment, the SUPG method modifies Eq. \eqref{1D-weak} to 
\begin{equation}
 \underbrace{(w^h,i\omega \phi^h) + (w^h,a\phi^h_{,x}) + (w^h_{,x},\kappa \phi^h_{,x})}_\text{Baseline Galerkin}  + \underbrace{\sum_e\Big(\tau a w^h_{,x},i\omega \phi^h + a\phi^h_{,x} - \kappa \phi^h_{,xx}\Big)_{\Omega_e}}_\text{The SUPG terms}  = 0,
\label{1D-SUPG}
\end{equation}
where the integrals under summation are performed in elements interior $\Omega_e$ given that $\phi^h_{,xx}$ is not defined on the element boundaries.
The parameter $\tau$ in Eq. \eqref{1D-SUPG} is formulated in a way such that the resulting approximate solution is nodally exact for steady-state flows with $\omega=0$.

\cola{For linear interpolation functions employed here, the diffusive term in the SUPG terms in Eq. \eqref{1D-SUPG} vanishes.
Nevertheless, this term is retained to ensure the formulation's adherence to its temporal counterpart.
Note that it is possible to estimate that diffusive term by performing $L_2$ projection of $\phi^h_{,x}$ to the mesh nodes and use the resulting variable to compute $\phi^h_{,xx}$ in the element interior. 
As detailed in Appendix \ref{app:CSUPG}, this modification generates less accurate results for the problem under consideration, hence it will not be considered in what follows.}

The two added SUPG terms effectively modify the convective velocity to $\tilde a = a - i\omega \tau a$ and the diffusion coefficient to $\tilde \kappa = \kappa + \tau a^2$. 
That effectively modifies $\alpha = ah/(2\kappa)$ to $\tilde \alpha = \tilde ah/(2\tilde \kappa)$. 
Therefore, the new numerical solution must be computed by substituting $\alpha$ with $\tilde \alpha$ in Eq. \eqref{rho12}. 
Recalling that $\beta =0$ in the steady regimes, that results in $\tilde \rho_{1,2}=1,(1+\tilde \alpha)/(1-\tilde \alpha)$. 
Thus, the SUPG solution computed based on $\tilde \rho_{1,2}$ will be nodally exact when~\cite{hughes1979multidimentional}  
\begin{equation}
    \tilde \rho_{1,2}^A = \exp(r_{1,2}\frac{x_A}{L}),
\label{SUPG-der}
\end{equation}
where $A$ is an exponent on the left-hand side. 

From the two conditions imposed by Eq. \eqref{SUPG-der}, one is already satisfied given that $r_1=0$ (since $W=0$) and $\tilde \rho_1=1$. 
It is from the second condition that one can obtain a relationship for $\tau$. 
In its exact form, that relationship is 
\begin{equation}
    \tau = \frac{h}{2a}\left(\coth\alpha - \frac{1}{\alpha}\right). 
\label{tau_e}
\end{equation}
To simplify computation and extension of Eq. \eqref{tau_e} to multiple dimensions, $\tau$ is often approximated as $\tau = \frac{h}{2a}(1+9\alpha^{-2})^{-\frac{1}{2}}$. 
The resulting $\tau$, which behaves the same as $\tau$ asymptotically, is often written as~\cite{hughes1986generalized,shakib1991new,bazilevs2007variational}
\begin{align}
    \tau &= (\tau_{\rm conv}^{-2} + \tau_{\rm diff}^{-2})^{-\frac{1}{2}}, \label{tau_supg}\\
    \tau_{\rm conv}^{-1} &= \frac{2a}{h}, \label{tau_conv}\\
    \tau_{\rm diff}^{-1} &= \frac{12\kappa}{h^2}. \label{tau_diff}
\end{align}

Although the above method is designed for the steady regime, one may directly apply it to the time-spectral form of the convection-diffusion equation. 
That entails computing $\tau$ from Eq. \eqref{tau_supg} and plugging it into Eq. \eqref{1D-SUPG} to compute the numerical solution $\phi^h$. 
Later in Section \ref{sec:1dtest}, we will discuss how this method behaves if it were to be used in unsteady regimes where $\omega \ne 0$. 

Lastly, we must point out that an alternative form of Eq. \eqref{tau_supg} has been proposed for space-time finite element methods to incorporate the effect of acceleration term in unsteady settings by adding a term that amounts to $\tau_{\rm acc}^{-2}=\omega^2$ for the problem under consideration \cite{shakib1991new}. 
The numerical experiments, which are not reported here, show the inferior performance of this design of $\tau$. 
Hence, in what follows, we forgo the discussion of this alternative form of Eq. \eqref{tau_supg}.

\subsection{Variational multiscale method}\label{sec:1d-vms}
The VMS method is constructed~\cite{hughes1995multiscale,hughes1998variational,hughes2007variational} by modeling the scale in $\phi$ that are not resolved by $\phi^h$, namely $\phi^\prime$, via the residual of the original PDE, i.e., $r(\phi^h)$. 
More specifically, to build this method we set
\begin{equation}
    \phi = \phi^h + \phi^\prime, 
\label{phi-vms}
\end{equation}
with 
\begin{equation}
    \phi^\prime = -\tau r(\phi^h), 
\label{phip}
\end{equation}
in which 
\begin{equation}
   r(\phi^h) = i\omega \phi^h + a \phi^h_{,x} - \kappa \phi^h_{,xx}.  
\label{res}
\end{equation}

With these definitions, the VMS problem statement becomes similar to that of the Galerkin's in Eq. \eqref{1D-weak} when $\phi^h$ is replaced with $\phi$ from Eq. \eqref{phi-vms}. 
The result is
\begin{equation}
 \underbrace{(w^h,i\omega \phi^h) + (w^h,a\phi^h_{,x}) + (w^h_{,x},\kappa \phi^h_{,x})}_\text{Baseline Galerkin}  + \underbrace{\sum_e\Big(\tau a w^h_{,x},r(\phi^h)\Big)_{\Omega_e}}_\text{The SUPG terms} - \underbrace{\sum_e\Big(i\omega \tau w^h,r(\phi^h)\Big)_{\Omega_e}}_\text{The new VMS terms} = 0. 
\label{1D-vms}
\end{equation}

One can make two observations by comparing the proposed time-spectral VMS method in Eq. \eqref{1D-vms} against its conventional temporal counterpart, which is discussed later in Section \ref{sec:patient}. The first is the similarity between the two where the diffusive flux of the unresolved scale is neglected in both. The second is the distinction between the two, where the acceleration of the unresolved scale is modeled only in the time-spectral VMS method. 
This contribution, which is captured by the new VMS terms in Eq. \eqref{1D-vms}, is entirely absent in the time formulation as $\hat u^\prime_{,t}$ is typically dropped from the discrete form. That simplification effectively reduces the conventional VMS method to the SUPG method. Later, we will see how the new VMS terms in Eq. \eqref{1D-vms} improve the accuracy and stability of this method.

\subsection{Galerkin/least-squares method}
The GLS method is constructed by incorporating a symmetric penalty term that is proportional to the residual of the original PDE to the baseline Galerkin's method \cite{hughes1989new,shakib1991new,shakib1989finite}. 
That penalty term is also scaled by $\tau$ to recover the steady SUPG term, thus producing
\begin{equation}
     \Big(w^h,r(\phi^h)\Big) + \sum_e\Big(r^*(w^h),\tau r(\phi^h)\Big)_{\Omega_e}= 0,
\label{lsq-base}
\end{equation}
where $f^*$ denotes the complex conjugate of $f$.

In constructing the penalty term in Eq. \eqref{lsq-base}, it is crucial to employ the adjoint operator $r^*(w^h)$. This way, the contribution of the penalty term to the tangent matrix will be a positive-definite Hermitian matrix. This property ensures the resulting method produces a well-posed linear system and a stable solution (Section \ref{sec:conv}).

Simplifying Eq. \eqref{lsq-base} for linear interpolation functions results in a method that is identical to the VMS method discussed in Section \ref{sec:1d-vms}. Thus, in what follows, we use the umbrella term VMS/GLS to label the results that are obtained from Eq. \eqref{1D-vms}.

\subsection{A new augmented SUPG method} \label{sec:asu}
Our overall approach to obtaining a new stabilization method, which we call augmented SUPG or ASU, for the time-spectral convection-diffusion equation is similar to that of the SUPG. 
The key difference is that, in this case, we enforce both conditions in Eq. \eqref{SUPG-der} while taking into account unsteady flows where $\beta\ne 0$. 

Enforcing nodally exact solution in general for unsteady flows entails computing two distinct $\hat \rho_{1,2}$ in terms of modified $\hat \alpha$ and $\hat \beta$ and making sure they satisfy a relationship similar to that of Eq. \eqref{SUPG-der}. 
The element Peclet and Womersley number, on the other hand, can be modified by properly adjusting oscillation frequency $\omega$, convective velocity $a$, and/or diffusivity $\kappa$. 
Of course, out of these three parameters, only two can be independently adjusted given the solution to the discrete form is insensitive up to a prefactor. 
Although the final result is independent of which two parameters we select, adjusting $\omega$ and $\kappa$ slightly simplifies the overall derivation process. 
Thus, the ASU method is derived by seeking $\hat \omega$ and $\hat \kappa$ or equivalently 
 \begin{align}
\hat \alpha &= \frac{ah}{2\hat \kappa},\label{alphah} \\
\hat \beta &= \frac{\hat \omega h^2}{6\hat \kappa}, 
\label{betah}
\end{align}
such that the discrete solution to the time-spectral convection-diffusion problem becomes nodally exact for piecewise linear shape functions. 

The process of computing $\hat \alpha$ and $\hat \beta$ so that the resulting $\hat \rho_{1,2}(\hat \alpha,\hat \beta)$ satisfy the corresponding conditions imposed by Eq. \eqref{SUPG-der} is rather lengthy, and thus, moved to Appendix \ref{app:pres}. 
The result of that process is
\begin{align}
\hat \alpha &= \frac{3\sinh\alpha}{\cosh\gamma + 2\cosh\alpha}, \label{alphah1}\\
i\hat \beta &= \frac{\cosh\gamma - \cosh\alpha}{\cosh\gamma + 2\cosh\alpha}, 
\label{betah1}
\end{align}
where 
\begin{equation}
\gamma = \sqrt{\alpha^2 + 6i\beta}. 
\label{gamma}
\end{equation}

One can readily verify that Eq. \eqref{alphah1} generalizes the corresponding relationship in the SUPG method as setting $\beta=0$ will produce the well-known relationship $\hat \alpha = \tanh\alpha$ while $\hat \beta = \beta$ is unchanged at zero. 

Knowing the desired forms of $\hat \alpha$ and $\hat \beta$, our next task is to design $\hat \omega$ and $\hat \kappa$ so that those desired forms are achieved. 
Since $\hat \kappa = \kappa \alpha /\hat \alpha$, from Eq. \eqref{alphah1} we can write  
\begin{equation}
\hat \kappa =  \kappa \alpha \coth\alpha + \kappa \alpha\left( \frac{\cosh \gamma  - \cosh\alpha }{3\sinh \alpha}\right). 
\label{kappah1}
\end{equation}
The first term on the right-hand side is identical to that of the traditional SUPG method and can be written as $\kappa + \tau a^2$ based on Eq. \eqref{tau_e}. 
To simplify the second term, note that from Eqs. \eqref{alphah1} and \eqref{betah1} 
\begin{equation}
\frac{\cosh \gamma  - \cosh\alpha }{3\sinh \alpha} = \frac{i\hat \beta}{\hat \alpha}. 
\label{kappah2}
\end{equation}
Relating this result to $\hat \omega$ and $\tau_{\rm diff}$ via Eqs. \eqref{alpha}, \eqref{tau_diff}, \eqref{alphah}, and \eqref{betah} produces 
\begin{equation}
\hat \kappa =\underbrace{\kappa}_\text{Baseline Galerkin} + \underbrace{a^2 \tau}_\text{Conventional SUPG} + \underbrace{\kappa_{_{\rm ASU}}}_\text{The new ASU term}  
\label{kappah}
\end{equation}
where 
\begin{equation}
\kappa_{_{\rm ASU}} = 2i \hat \omega \tau_{\rm diff} \kappa. 
\label{kasu}
\end{equation}
This strikingly simple expression suggests that although the SUPG method remains mostly intact for the time-spectral form of the convection-diffusion equation, its diffusion must be augmented by $\kappa_{_{\rm ASU}}$. 
Note that $\kappa_{_{\rm ASU}}$ is imaginary to the leading order as it is pre-multiplied by $i$ and, as we will see later, $\hat \omega$ is real up to the first order with respect to $\beta$. 
That implies $\kappa_{_{\rm ASU}}$ primarily acts to diffuse the scalar field between its real and imaginary components (i.e., in-phase and out-of-phase solutions).

Similar to $\hat \kappa$, one can compute $\hat \omega$ by observing $\hat \omega/\omega = (\hat \beta \alpha)/(\beta \hat \alpha)$.
That yields 
\begin{equation}
\hat \omega = \left(\frac{\alpha}{i\beta}\right) \left(\frac{\cosh \gamma  - \cosh\alpha }{3\sinh \alpha}\right) \omega. 
\label{omegahe}
\end{equation}
This seemingly complex relationship can be significantly simplified if $\beta \lessapprox 1$. 
As detailed in Appendix \ref{omegah-app}, using asymptotic expansions at the limits $\beta/\alpha^2 \ll 1$ and $\beta/\alpha^2 \gg 1$ produces another strikingly simple relationship for $\hat \omega$, that is
\begin{equation}
\hat \omega \approx \omega \exp(i\omega \tau ),
\label{omegah}
\end{equation}
in which $\tau$ is the approximate SUPG time scale from Eq. \eqref{tau_supg} that is already calculated when implementing any of the stabilized methods discussed above. 

\begin{figure}
  \centering
  \include{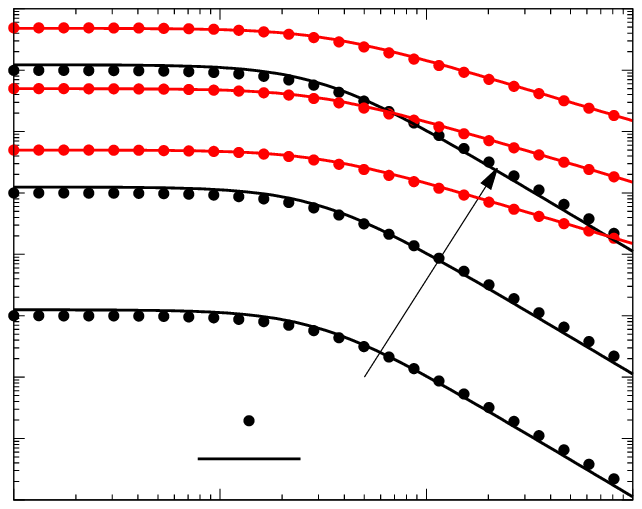}
   \caption{The real (black) and imaginary (red) component of $1-\hat\omega/\omega$ using the exact (dots) and approximate (solid line) expressions as a function of $\alpha$ for $\beta=0.01$, 0.1, and 1.0. Both real and imaginary components monotonically increase as $\beta$ increases.}
  \label{fig:tauo}    
\end{figure}

The leading order term in the relationship for $\hat \omega$ (Eq. \eqref{omegah}) with regard to $\beta$ (viz., a measure of flow unsteadiness relative to diffusion at element scale) is $\omega$.
In other words, the exponential prefactor in Eq. \eqref{omegah} goes to one as $\omega\to 0$ or $\beta\to0$.
This observation is, of course, expected by design as $\omega$ must not be altered at the steady state limit. 
Thus, to highlight minute differences between the exact and approximate $\hat \omega$ from Eqs. \eqref{omegahe} and \eqref{omegah}, respectively, we subtracted this leading order term from the two before comparing the two in Figure \ref{fig:tauo}.
Although shown only up to $\beta =1$, the approximate $\hat \omega$ closely agrees with its exact form for $\beta \lessapprox 5$. 
For $\beta$ values higher than that, Eq. \eqref{omegahe} exhibits a highly nonlinear behavior that is not captured by Eq. \eqref{omegah}. 
However, finding an alternative approximate form for $\hat \omega$ that captures this extreme regime is of little interest as it signifies situations where the mesh resolution is extremely poor in comparison to the level of details present in the solution. 

As it is evident in Figure \ref{fig:tauo}, $\hat \omega$, passed its leading order term, follows a behavior similar to that of $\tau$. 
This well-known behavior is characterized by the diffusive and convective limits at which $\tau \approx \tau_{\rm diff}$ and $\tau \approx \tau_{\rm conv}$, respectively. 
Similarly, in the diffusive limit at which $\alpha \ll 1$, $\hat \omega$ becomes independent of $\alpha$. 
In this limit, $\hat \omega/\omega$ passed its constant leading order becomes proportional to $i\beta$ (Eq. \eqref{C6}).
As we will see in the next section, it is at this limit that the ASU provides a more accurate estimate than the conventional SUPG method. 
The two methods, however, both converge in the convective limit at which $\alpha \gg 1$ since the difference between $\omega$ and $\hat \omega$ diminishes at a rate proportional to $\alpha^{-1}$ (Eq. \eqref{C3}). 

Having a relationship for $\hat \omega$ and $\hat \kappa$, the ASU method is obtained by modifying the baseline Galerkin's method to 
\begin{equation}
 \underbrace{(w^h,i\hat \omega \phi^h) + (w^h,a\phi^h_{,x}) + (w^h_{,x},\kappa \phi^h_{,x})}_\text{Galerkin with a modified $\omega$}  + \underbrace{(\tau a w^h_{,x},a\phi^h_{,x})}_\text{The steady SUPG term} + \underbrace{(w^h_{,x},\kappa_{_{\rm ASU}}\phi^h_{,x})}_\text{The new ASU term} = 0, 
\label{1D-ASU}
\end{equation}
in which $\hat \omega$ and $\kappa_{_{\rm ASU}}$ are computed from Eqs. \eqref{omegah} and \eqref{kasu}, respectively. 

Considering Eq. \eqref{1D-ASU}, the simplest way to implement the ASU method is to  
\begin{enumerate}
\item start from the Galerkin's method but use $\hat \omega$ instead of $\omega$ throughout computations, 
\item add the quasi-steady SUPG term, and 
\item augment the physical diffusivity with the dominantly-imaginary $\kappa_{_{\rm ASU}}$. 
\end{enumerate}

In the next section, we will see how the five methods discussed above behave in the 1D unsteady convection-diffusion setting. 
Later in Section \ref{sec:md}, we will discuss how they can be generalized to multiple dimensions and behave in more complex settings.  

\subsection{Comparison of various methods in 1D} \label{sec:1dtest}
All the methods discussed above can be brought into a form similar to that of Galerkin's method in Eq. \eqref{1D-weak} but with a potentially modified $\omega$, $a$, and $\kappa$. 
Once they are in this standard form,  we can identify how each alters the baseline frequency, velocity, and diffusivity parameters as they appear in Galerkin's discrete form.  
The result of this process is condensed in Table \ref{table:sum}.
The following are the key observations: 

\begin{enumerate}
    \item The quasi-steady SUPG term appears in all stabilized methods, thus increasing the physical diffusivity by $a^2\tau$. 
    The inclusion of this term is crucial in obtaining stable results at highly convective regimes. 
    This observation is in accordance with what has been traditionally observed in the time formulation of convection-diffusion equation~\cite{hughes1986generalized}.
    \item All stabilization methods will be similar in strongly convective regimes that are not highly oscillatory ($\alpha\gg1$ and $\beta \lessapprox 1$). 
    In such regimes, the overall behavior of these methods is determined by $a^2\tau$ term in $\hat \kappa$ rather than $i\omega\tau$.
    \item All methods agree on how $\omega$ must be adjusted up to the second order leading term with respect to $\beta$ (or $\omega \tau$).
    If we divide the SUPG and VMS/GLS discrete forms by $1-i\omega\tau$ and $1-2i\omega\tau$, respectively, so that all methods leave convective velocity unchanged at $\hat a= a$, we can conclude that $\hat \omega/\omega = 1+i\omega \tau + O\Big((\omega\tau)^2\Big)$ for all stabilization methods. 
    By the same virtue, the VMS/GLS and ASU agree on the form of the effective diffusivity $\hat \kappa$ at small $\alpha$ when $a^2 \tau$ is negligible (as stated earlier, at large $\alpha$ all stabilization methods collapse). This observation explains the lower accuracy of the SUPG method in comparison to the VMS/GLS and ASU, as shown below.
\end{enumerate}

\begin{table}
\begin{center}
\begin{tabular}{c|cccc}
  & Galerkin & SUPG & VMS/GLS & ASU  \\  \hline\hline
    $\hat \omega$ & $\omega$ & $\omega$           & $(1-i\omega\tau)\omega$  & $\exp(i\omega\tau)\omega$                        \\
    $\hat a$      & $a$      & $(1-i\omega\tau)a$ & $(1-2i\omega\tau)a$      & $a$                                              \\ 
    $\hat \kappa$ & $\kappa$ & $\kappa + a^2\tau$ & $\kappa+a^2\tau$         & $\kappa+2i\hat \omega \tau_{\rm diff} + a^2\tau$ \\\hline \hline
\end{tabular}
\caption{The effective oscillation frequency $\hat \omega$, convective velocity $\hat a$, and diffusivity $\hat \kappa$ for a given method if they were to be formulated in the form of Eq. \eqref{1D-weak}.}
\label{table:sum}
\end{center}
\end{table} 

\begin{figure}
  \centering
  \include{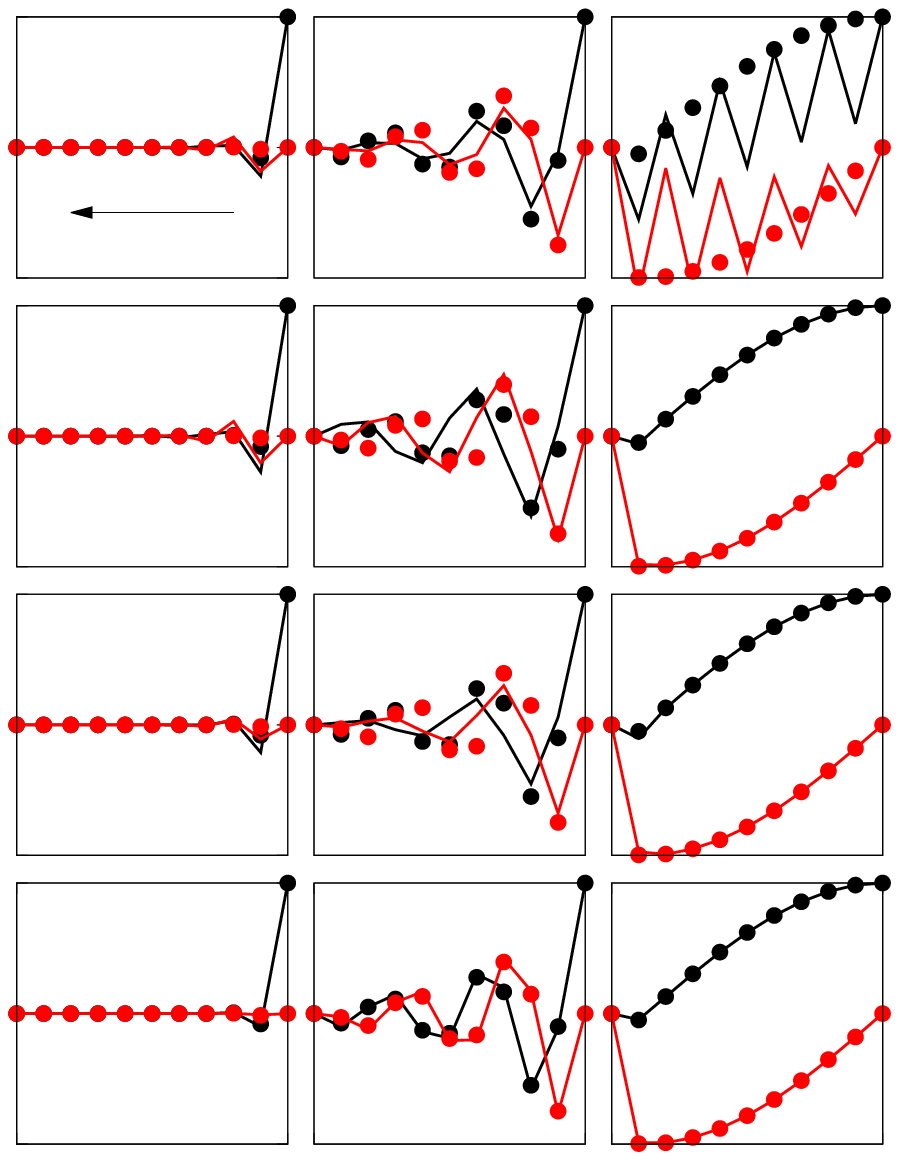}
   \caption{The real (black) and imaginary (red) components of the solution for the 1D model problem using various methods (solid line) in comparison to the exact solution (dots).  $\beta=3$ for all cases and $\alpha$ is 0.5, 5, and 50 for the left, middle, and right columns, respectively. }
  \label{fig:1d}
\end{figure}

\begin{figure}[H]
  \centering
  \include{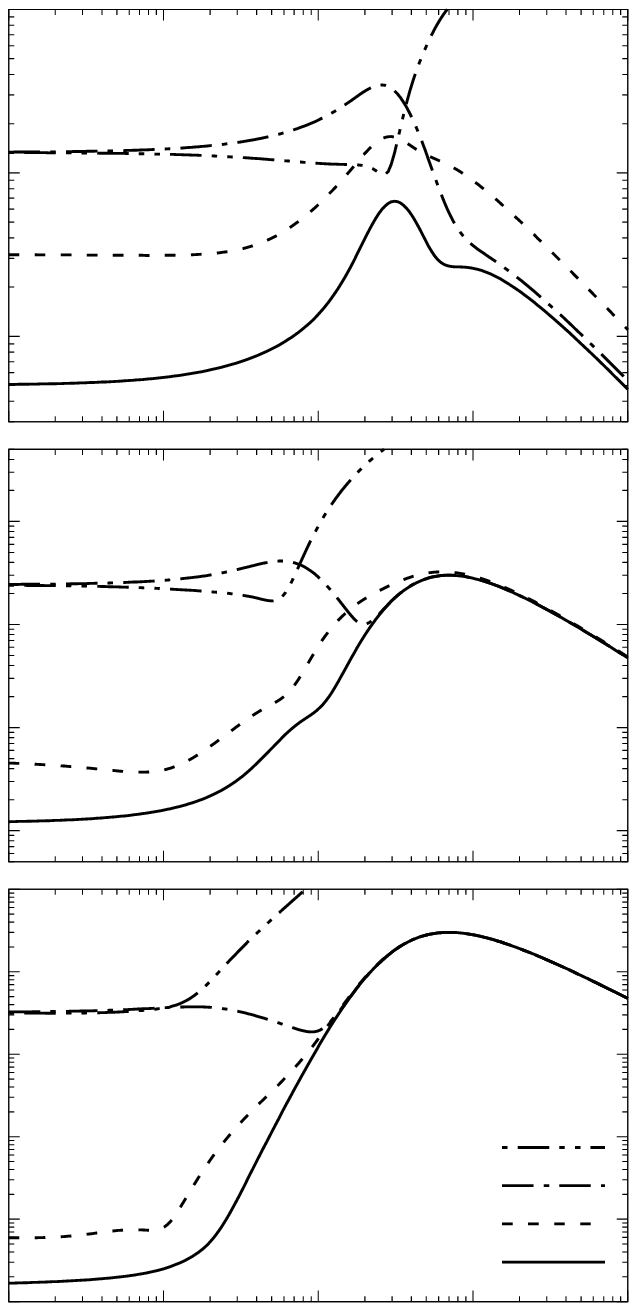}
   \caption{The numerical solution error for the 1D model problem as a function of $\alpha$ for (a) $\beta=1$, (b) $\beta=0.1$, and (c) $\beta=0.01$ using the Galerkin's (dash-double-dotted), SUPG (dash-dotted), VMS/GLS (dotted), and ASU (solid) approach. }
  \label{fig:err}    
\end{figure}

To put the above observations into a more concrete perspective, we have simulated the 1D model problem described earlier in Section \ref{sec:1D-problem-statement} using all five methods. 
The results are shown in Figure \ref{fig:1d} for three isolated cases at which $\beta=3$ and $\alpha=0.5$, 5, and 50. 
To evaluate the accuracy of these methods more comprehensively, we have repeated these computations for many values of $\alpha$ and calculated the $L_2$-norm of error for each simulation. 
The results are reported in Figure \ref{fig:err} for $\beta=0.01$, 0.1, and 1. 

The general observations from these numerical results agree with what we discussed above. Namely 
\begin{enumerate}
    \item All stabilization methods (i.e., SUPG, VMS/GLS, and ASU) will be similar at sufficiently large $\alpha$, producing a reasonably accurate result. The Galerkin's method, as we saw earlier in Section \ref{sec:issue}, suffers from nonphysical oscillations in such regimes. 
    \item For sufficiently small $\beta$, the VMS/GLS and ASU methods behave similarly. 
    Additionally, the Galerkin and SUPG behave similarly if both $\alpha$ and $\beta$ are small. In such regimes, the first group generates more accurate results than the second group. 
    \item Overall, the ASU is the most accurate method as it is tailor-designed for this problem. 
    The error associated with this method is purely due to the approximations associated with Eqs. \eqref{tau_supg} and \eqref{omegah}. 
    Differently put, the error shown in Figure \ref{fig:err} for the ASU method reduces to zero if $\tau$ and $\hat \omega$ are calculated exactly from Eqs. \eqref{tau_e} and \eqref{omegahe}.
    Our numerical experiments show that the majority of the error is due to the approximate $\tau$ (Eq. \eqref{tau_supg}). 
    It is only at relatively high $\beta$ that the approximation used in calculating $\hat \omega$ (Eq. \eqref{omegah}) translates to an error in $\phi^h$.
    \item The second most accurate method, after the ASU, varies depending on the regime under consideration. Roughly speaking, the VMS/GLS produces accurate results except for large $\beta$ and small $\alpha$. 
    That is followed by the SUPG which remains reasonably accurate for the bulk of the parameter space. 
    The same can be said for the Galerkin's method if we exclude high $\alpha$. 
\end{enumerate}

\section{Generalization to multiple dimensions} \label{sec:md}
In what follows, we briefly state the discrete form of the methods considered in the previous section and discuss their practical implementation in a finite element code that uses purely real arithmetic. 
We do so since the majority of existing codes use purely real variables, and they are linked against linear solvers that are purposely built for handling real arithmetic. 

For modeling the transport of a non-reactive solute in multiple dimensions, we consider the general form of Eq. \eqref{1D-fourier} that is 
\begin{equation}
  \begin{alignedat}{3}
    i\omega \phi + \bl a \cdot \nabla \phi &= \nabla \cdot (\kappa \nabla \phi) + q \;\;\; & &   {\rm in } \;\;\; & & \Omega,  \\ 
    \phi &= g & & {\rm on} & & \Gamma_g,  \\ 
    \kappa \bl n \cdot \nabla \phi &= h & & {\rm on} & & \Gamma_h, 
\end{alignedat}
\label{ad_md}
\end{equation}
where $\bl n$ is the boundary outward normal vector, $q(\bl x)\in \mathbb Z$ is the given source term, and $g(\bl x)\in \mathbb Z$ and $h(\bl x)\in \mathbb Z$ are the prescribed Dirichlet and Neumann boundary conditions, respectively. 

When formulating Eq. \eqref{ad_md}, we did not incorporate a $s\phi$ reaction term on the right-hand side that scales with the solution. 
Nonetheless, by considering Eq. \eqref{ad_md}, we are also considering such scenarios as they are also governed by Eq. \eqref{ad_md} when $i\omega$ is replaced by $i\omega - s$. 
In fact, the ASU method introduced in this study can be adopted as a new stabilization method for the convection-diffusion-reaction equation that is formulated in time. 
The resulting method, which is obtained by replacing $i\omega$ with $-s$ in all relevant expressions, will also exhibit a super-convergence behavior by generating a nodally exact solution for the 1D model problem of Section \ref{sec:1d}. 
Consistent with earlier observations, our numerical experiments involving time simulation of the convection-diffusion-reaction equation show that the ASU method is more accurate than all other stabilization methods that we have considered in this study. 
Therefore, much of what is discussed here can be simply generalized to the time-spectral form of the convection-diffusion-reaction equation. 

To obtain a formulation that solely relies on real arithmetic, we reformulate Eq. \eqref{ad_md} in terms of its real and imaginary components. 
That requires taking the discrete solution as 
\begin{equation}
\phi^h = \phi^h_r + i\phi^h_i,
\label{chngov}
\end{equation}
plugging it in the discrete form, separating real and imaginary elements, and testing each with separate test functions.
The same decomposition is also applied to the source term and boundary conditions by taking 
\begin{equation}
\begin{split}
    q &=q_r+iq_i, \\
    g &= g_r + ig_i, \\
    h &= h_r + ih_i, \\  
\end{split}  
\label{qgh_decom}
\end{equation}
so that $q_r(\bl x)$, $q_i(\bl x)$, $g_r(\bl x)$, $g_i(\bl x)$, $h_r(\bl x)$, and $h_i(\bl x)$ are all real functions. 
Using these definitions, the baseline Galerkin's method can be stated as finding $\phi^h_r$ and $\phi^h_i$ that are equal to $g_r$ and $g_i$, respectively, on $\Gamma_g$ such for any  $w_r^h$ and $w_i^h$ that vanish on $\Gamma_g$, we have 
\begin{equation}
    B_{\rm G}(w^h_r, w^h_i; \phi^h_r, \phi^h_i) = F(w^h_r,w^h_i),
\label{gal-bg}
\end{equation}
where
\begin{equation}
\begin{split}
    B_{\rm G}(w^h_r, w^h_i; \phi^h_r, \phi^h_i) = &-(w^h_r,\omega \phi^h_i) +(w^h_r,\bl a\cdot  \nabla \phi^h_r) + (\nabla w^h_r,\kappa \nabla \phi^h_r) \\
    &+ (w^h_i,\omega \phi^h_r ) + (w^h_i,\bl a\cdot \nabla \phi^h_i) + (\nabla w^h_i,\kappa \nabla \phi^h_i). 
\end{split}
\label{gal}
\end{equation}
and 

\begin{equation}
    F(w^h_r, w^h_i) = (w^h_r,q_r) +(w^h_i,q_i) + (w^h_r,h_r)_{\Gamma_h} + (w^h_i, h_i)_{\Gamma_h}. 
\label{gal-RHS}
\end{equation}

Similarly, the SUPG method translates to finding $\phi^h_r$ and $\phi^h_i$ such for any  $w_r^h$ and $w_i^h$ we have 
\begin{equation}
    B_{\rm S}(w^h_r, w^h_i; \phi^h_r, \phi^h_i) = F(w^h_r,w^h_i).
\label{supg-bs}
\end{equation}
where
\begin{equation}
    B_{\rm S}(w^h_r, w^h_i; \phi^h_r, \phi^h_i) = B_{\rm G} + \sum_e \left[\Big(\tau \bl a\cdot \nabla w^h_r,r_r(\phi^h)\Big)_{\Omega_e} + \Big(\tau \bl a\cdot \nabla w^h_i,r_i(\phi^h)\Big)_{\Omega_e}\right],
\label{supg}
\end{equation}
and
\begin{equation}
\begin{split}
    r_r(\phi^h) &= -\omega \phi^h_i+\bl a\cdot  \nabla \phi^h_r - \nabla \cdot(\kappa \nabla  \phi^h_r) - q_r, \\
    r_i(\phi^h) &= +\omega \phi^h_r+\bl a\cdot  \nabla \phi^h_i - \nabla \cdot(\kappa \nabla  \phi^h_i) - q_i. \\
\end{split}
\label{r_ri}
\end{equation}

The stabilization parameter $\tau$ in Eq. \eqref{supg} is computed using a relationship similar to its 1D counterpart from Eq. \eqref{tau_supg} that is \cite{bazilevs2013computational,shakib1991new,bazilevs2007variational} 
\begin{align}
    \tau &= (\tau_{\rm conv}^{-2} + \tau_{\rm diff}^{-2})^{-\frac{1}{2}}, \label{tau_supg-md} \\
    \tau_{\rm conv}^{-2} &= \bl a^T \bl G \bl a, \label{tau_conv-md} \\
    \tau_{\rm diff}^{-2} &= 9\kappa^2 \bl G : \bl G, \label{tau_diff-md} 
\end{align}
where 
\begin{equation}
    \bl G = \left(\frac{\partial \bl \xi}{\partial \bl x}\right)^T\left(\frac{\partial \bl \xi}{\partial \bl x}\right), 
    \label{G_def}
\end{equation}
is the metric tensor computed from mapping between physical $\bl x$ and parent element $\bl \xi$ coordinates. 
The coefficient 9 in Eq. \eqref{tau_diff-md}, which may be optimized for a given element type, is set to that constant here to reduce the number of variables that could influence the relative accuracy of various methods. 

The generalization of the VMS/GLS method from Eq. \eqref{1D-vms} to multiple dimensions is also straightforward and amounts to finding $\phi^h_r$ and $\phi^h_i$ such for any  $w_r^h$ and $w_i^h$ we have 
\begin{equation}
    B_{\rm V}(w^h_r, w^h_i; \phi^h_r, \phi^h_i)=F(w^h_r,w^h_i).
\label{vms-bv}
\end{equation}
where 
\begin{equation}
    B_{\rm V}(w^h_r, w^h_i; \phi^h_r, \phi^h_i) = B_{\rm S}  + \sum_e \left[ \Big(\tau\omega w^h_r,  r_i(\phi^h) \Big)_{\Omega_e} - \Big(\tau \omega w^h_i,r_r(\phi^h)\Big)_{\Omega_e} \right].
\label{vms}
\end{equation}

Finally, the generalization of the ASU from Eq. \eqref{1D-ASU} is stated as finding $\phi^h_r$ and $\phi^h_i$ such for any $w_r^h$ and $w_i^h$ we have 
\begin{equation}
    B_{\rm A}(w^h_r, w^h_i; \phi^h_r, \phi^h_i)=F(w^h_r,w^h_i).
\label{asu-ba}
\end{equation}
where 
\begin{equation}
\begin{split}
    B_{\rm A}(w^h_r, w^h_i; \phi^h_r, \phi^h_i) = \hat B_{\rm G} &+(\tau \bl a\cdot \nabla w^h_r,\bl a \cdot\nabla \phi^h_r) + (\nabla w^h_r,\kappa_{_{\rm ASU}} \nabla \phi^h_r) \\
    &+(\tau \bl a\cdot \nabla w^h_i,\bl a\cdot\nabla \phi^h_i) + (\nabla w^h_i,\kappa_{_{\rm ASU}} \nabla \phi^h_i).
\end{split}
\label{asu1}
\end{equation}
In this equation, $\hat B_{\rm G}$ is identical to $B_{\rm G}$ except for $\omega$ that is replaced by $\hat \omega$ from Eq. \eqref{omegah} that is $\hat \omega = \exp(i\omega\tau)\omega$. 
Also, $\kappa_{_{\rm ASU}}$ is computed from Eq. \eqref{kasu} with $\tau_{\rm diff}$ from Eq. \eqref{tau_diff-md}, yielding
\begin{equation}
\kappa_{_{\rm ASU}} = \frac{2i }{3}\left(\bl G:\bl G\right)^{-\frac{1}{2}}\hat \omega. 
\label{kasu-md}
\end{equation}

Since $\hat \omega$ and $\kappa_{_{\rm ASU}}$ are complex-valued variables, Eq. \eqref{asu1} must be further simplified to obtain a purely real expression. 
Denoting the real and imaginary component of $\hat \omega$ by $\hat \omega_r$ and $\hat \omega_i$, respectively, and those of $\kappa_{_{\rm ASU}}$ by $\kappa_{r_{\rm ASU}}$ and $\kappa_{i_{\rm ASU}}$, respectively, the final expression for the ASU method becomes
\begin{equation}
\begin{split}
    B_{\rm A}(w^h_r, w^h_i; \phi^h_r, \phi^h_i) = &-(w^h_r,\hat \omega_r \phi^h_i)-(w^h_r,\hat \omega_i \phi^h_r) +(w^h_r,\bl a\cdot  \nabla \phi^h_r) + (\nabla w^h_r,\kappa \nabla \phi^h_r) \\
    &+ (w^h_i,\hat \omega_r \phi^h_r ) -(w^h_i,\hat \omega_i \phi^h_i) + (w^h_i,\bl a\cdot \nabla \phi^h_i) + (\nabla w^h_i,\kappa \nabla \phi^h_i)  \\
    &+(\tau \bl a\cdot \nabla w^h_r,\bl a \cdot\nabla \phi^h_r) + (\nabla w^h_r,\kappa_{r_{\rm ASU}} \nabla \phi^h_r) - (\nabla w^h_r,\kappa_{i_{\rm ASU}} \nabla \phi^h_i) \\
    &+(\tau \bl a\cdot \nabla w^h_i,\bl a\cdot\nabla \phi^h_i) + (\nabla w^h_i,\kappa_{r_{\rm ASU}} \nabla \phi^h_i)+ (\nabla w^h_i,\kappa_{i_{\rm ASU}} \nabla \phi^h_r).
\end{split}
\label{asu}
\end{equation}


\subsection{Stability properties} \label{sec:conv}
Before we test these methods in more realistic 2D and 3D settings, it will be instrumental to investigate their convergence behavior.
For a method to be convergent, it must be consistent and stable. 
The consistency of various methods is demonstrated in Appendix \ref{app:cons} and later in Section \ref{sec:msh}. 
Thus, for the remainder of this section, we focus on investigating the stability of these methods as the main requirement for their convergence.

Provided that the problem under consideration is a boundary value problem, we investigate the stability of various methods in terms of the properties of their tangent matrix~\cite{reinhardt2012analysis}. 
More specifically, we will examine whether the linear system produced by these methods is non-singular or avert to ill-conditioning.
Such a property has implications on the ability of the underlying iterative linear solver to produce a unique and stable solution using relatively few iterations. 

To ensure the stability of a method, one may study the positive definiteness of the underlying tangent matrix. 
For the baseline Galerkin's method, the tangent matrix will be positive definite if  
\begin{equation} 
    E_{\rm G} = \bl c^{\rm T}\bl K_{\rm G} \bl c, 
    \label{pd-req}
\end{equation}
is always larger or equal to zero for any $\bl c \in \mathbb R^N$ and it is zero if and only if $\bl c = 0$.
In Eq. \eqref{pd-req}, $\bl K_{\rm G}$ denotes the tangent matrix obtained from the Galerkin's method. 
It is rather straightforward to show that 
\begin{equation} 
    E_{\rm G} = B_{\rm G}(w^h_r,w^h_i;w^h_r,w^h_i), 
    \label{pd-req2}
\end{equation}
given that $\bl K_{\rm G}$ is extracted from Eq. \eqref{gal} and $w^h_r$ and $w^h_i$ are also arbitrary test functions, enabling us to evaluate Eq. \eqref{pd-req} for an arbitrary $\bl c$.  
Therefore, from Eqs. \eqref{gal} and \eqref{pd-req2} we have 
\begin{equation}
\begin{split}
    E_{\rm G} &= -(w^h_r,\omega w^h_i) +(w^h_r,\bl a\cdot  \nabla w^h_r) + (\nabla w^h_r,\kappa \nabla w^h_r) + (w^h_i,\omega w^h_r ) + (w^h_i,\bl a\cdot \nabla w^h_i) + (\nabla w^h_i,\kappa \nabla w^h_i) \\   
    &=  \int_{\Omega} \kappa  \|\nabla w^h\|^2 {\rm d}\Omega \ge 0,  
\end{split}
\label{cv-gal}
\end{equation}
where $w^h = [w_r^h\; w^h_i]^{\rm T}$.  
Since $E_{\rm G} = 0 $ only when $w^h=0$, the tangent matrix produced by the baseline Galerkin's method is formally positive definite. 
In practice, however, the resulting system can become ill-conditioned at high frequencies or in convection-dominated regimes since $E_{\rm G}$ does not depend on $\omega$ or $\bl a$, respectively. 
Thus, we expect the Galerkin's method to remain stable but produce an ill-conditioned matrix when $\beta\gg 1$ or $\alpha \gg 1$.
Our numerical experiment results corroborate this conclusion. 

The same procedure can be applied to the remaining stabilized methods. 
For the SUPG method, it follows from Eq. \eqref{supg} for linear interpolation functions
\begin{equation}
\begin{split}
    E_{\rm S} & = E_{\rm G}  + (\tau \bl a\cdot \nabla w^h_r,\bl a\cdot  \nabla w^h_r-\omega w^h_i) + (\tau \bl a\cdot \nabla w^h_i,\bl a\cdot  \nabla w^h_i+\omega w^h_r) \\
    &= \int_{\Omega} \Big[ \kappa  \|\nabla w^h\|^2 + \tau  \|\bl a\cdot\nabla w^h\|^2 + 2\tau\omega  w^h_r \bl a\cdot \nabla w^h_i \Big] {\rm d}\Omega. 
\end{split}
\label{cv-supg}
\end{equation}
While the first and second terms under the integral in Eq. \eqref{cv-supg} are always positive for a nonzero $w^h$, the third term could be positive or negative. 
$E_{\rm S}$ can, in fact, become negative if $\omega > O(a/h)$ or $\omega > O(\kappa^2/(ah^3))$ in strongly convective or diffusive regimes, respectively, thus introducing negative eigenvalues in the tangent matrix. 
With some eigenvalues being positive and some being negative, an eigenvalue may land near zero (given the significant variability in the mesh size, for instance), thus preventing a stable solution.   

Starting from Eq. \eqref{vms} and applying the same analysis to the VMS/GLS method yields 
\begin{equation}
\begin{split}
    E_{\rm V} = E_{\rm G} &+ (\tau \omega w^h_r,\omega w^h_r) + (\tau \bl a\cdot \nabla w^h_r,\bl a\cdot  \nabla w^h_r) - (2\tau \bl a\cdot \nabla w^h_r,\omega w^h_i)\\
    &+(\tau \omega w^h_i,\omega w^h_i) + (\tau \bl a\cdot \nabla w^h_i,\bl a\cdot  \nabla w^h_i) + (2\tau \bl a\cdot \nabla w^h_i,\omega w^h_r) \\
     = E_{\rm G} & + \int_{\Omega} \tau\Big[\|\omega w^h_i\|^2 + \|\bl a\cdot\nabla w^h_r\|^2 - 2\omega w^h_i \bl a\cdot \nabla w^h_r \Big] {\rm d}\Omega\\
     & +\int_{\Omega} \tau\Big[\|\omega w^h_r\|^2 + \|\bl a\cdot\nabla w^h_i\|^2 + 2\omega w^h_r\bl a\cdot \nabla w^h_i\Big] {\rm d}\Omega \\ 
     =\int_{\Omega} &\Big[ \kappa  \|\nabla w^h\|^2 + \tau\|\omega w^h_i - \bl a\cdot \nabla w^h_r\|^2 + \tau \|\omega w^h_r + \bl a\cdot\nabla w^h_i\|^2 \Big] {\rm d}\Omega \ge 0.   
\end{split}
\label{cv-vms}
\end{equation}
Therefore, the VMS/GLS method, similar to Galerkin's method, formally produces a positive definite tangent matrix and is expected to remain stable. 
That said, the VMS/GLS resulting tangent matrix is expected to be better conditioned at high $\alpha$ and $\beta$ regimes since $E_{\rm V}$ will scale with $\tau \|\bl a\cdot \nabla w^h\|^2$ and $\tau\|\omega w^h\|^2$, respectively, to remain sufficiently large.

Lastly, we use Eq. \eqref{asu} to evaluate the ASU.
The result is 
\begin{equation}
\begin{split}
E_{\rm A} &= -(w^h_r,\hat \omega_r w^h_i) +(w^h_r,\hat \omega_i w^h_r) + (w^h_r,\bl a\cdot  \nabla w^h_r) + (\nabla w^h_r,\kappa \nabla w^h_r) \\
     &\;\;\;\; + (w^h_i,\hat \omega_r w^h_r )+ (w^h_i,\hat \omega_i w^h_i) + (w^h_i,\bl a\cdot \nabla w^h_i) + (\nabla w^h_i,\kappa \nabla w^h_i) \\
   &\;\;\;\;+(\tau \bl a\cdot \nabla w^h_r,\bl a \cdot\nabla w^h_r) + (\nabla w^h_r,\kappa_{r_{\rm ASU}} \nabla w^h_r)  - (\nabla w^h_r,\kappa_{i_{\rm ASU}} \nabla w^h_i) \\
    &\;\;\;\;+(\tau \bl a\cdot \nabla w^h_i,\bl a\cdot\nabla w^h_i) + (\nabla w^h_i,\kappa_{r_{\rm ASU}} \nabla w^h_i)+ (\nabla w^h_i,\kappa_{i_{\rm ASU}} \nabla w^h_r) \\
    &= \int_{\Omega} \Big[\hat \omega_i  \|w^h\|^2 + \kappa \|\nabla w^h\|^2 +  \tau \|\bl a\cdot \nabla w^h\|^2 + \kappa_{r_{\rm ASU}} \|\nabla w^h\|^2\Big] {\rm d}\Omega. 
\end{split}
\label{cv-asu}
\end{equation}
At the outset, the ASU may appear to produce a positive definite tangent matrix as well. 
Nonetheless, that is not necessarily the case as $\hat \omega_i$ and $\kappa_{r_{\rm ASU}}$ may become negative. 
While the former only occurs at $\omega\tau > \pi$, the latter occurs even at small values of $\omega$. 
However, for the system to potentially produce negative eigenvalues in those scenarios, we must have $\kappa + \kappa_{r_{\rm ASU}} < 0$, which implies $\tau \omega^2 > O(\kappa/h^2)$. 
Thus, in either case, one expects the ASU to produce a positive definite tangent matrix at small $\omega$. 
We can, in fact, show that the ASU is limited by $\omega < O(a/h)$ and $\omega < O(\kappa/h^2)$ in strongly convective and diffusive regimes, respectively, to remain stable. 

To lessen those requirements, one can limit the angle $\omega\tau$ that appears in the definition of $\hat \omega$ in Eq. \eqref{omegah}. 
One way to achieve that will be to set an upper bound on $\tau$, namely $\tau_{\rm max}$, and compute $\hat \omega$ using the following expression rather than Eq. \eqref{omegah}:
\begin{equation}
    \hat \omega = \omega \exp(i\omega \min(\tau, \tau_{\rm max})).  
    \label{omegah-lim}
\end{equation}

Our numerical experiments show that $\tau_{\rm max}$ is set by the diffusive rather than the convective limit. 
Thus, considering the first two terms under the integral in Eq. \eqref{cv-asu}, we must have $|\hat \omega_i|<O(\kappa/h^2)$ for the ASU to remain stable. 
Since the first order approximation of $\hat \omega_i$ is $\omega^2 \tau$, the above condition translates to $\tau <O(\kappa/(h^2\omega^2))$, which provides us with an upper limit on $\tau$, namely $\tau_{\rm max}$. 
Using Eq. \eqref{tau_diff} to express $\kappa/h^2$ in terms of $\tau_{\rm diff}$ produces
\begin{equation}
\tau_{\rm max}^{-1} = \pi \omega^2 \tau_{\rm diff},  
    \label{tau_max}
\end{equation}
where $\tau_{\rm diff}$ is computed from Eq. \eqref{tau_diff-md} in a multi-dimensional setting. 
The pre-factor $\pi$ incorporated in Eq. \eqref{tau_max} is purely empirical and selected to be the smallest value that ensures stability of the ASU for a wide range of $\alpha$ and $\beta$. 
It is important to note that this value is based on numerical calculations using tetrahedral elements and may differ if one were to use a different interpolation function.

In summary, for the methods under consideration, only the baseline Galerkin's and VMS are formally guaranteed to produce a positive definite tangent matrix. 
The remaining methods, namely the SUPG and ASU could become unstable in strongly oscillatory regimes with a large Womersley number. 
Whether limiting the value of $\hat \omega_i$ in an ad hoc manner ensures the stability of the ASU is what we investigate below through a 2D and 3D test case. 
These cases also permit us to evaluate the accuracy of these methods in more realistic settings. 

\subsection{A 2D test case}
In Section \ref{sec:asu}, we discussed how one should modify the oscillations frequency $\omega$, convective velocity $a$, and diffusivity $\kappa$ to obtain a nodally exact solution in a 1D setting. 
In 2D and 3D settings, these modifications may not be ideal as they do not account for the direction-dependency of these parameters. 
Considering the conventional SUPG method, for instance, it is constructed to increase diffusion in the streamwise direction, leaving the diffusion in the crosswind direction unchanged. 
This selective modification of $\kappa$ is crucial as an increase in $\kappa$ in the crosswind direction will overly dampen the solution. 

\begin{figure}
  \centering
  \includegraphics[width=0.3\textwidth]{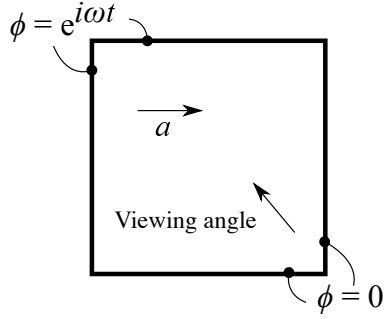}
   \caption{The 2D convection-diffusion problem under consideration in a square-shaped domain of size $L\times L$}
  \label{fig:2sch}
\end{figure}

To amplify these effects, we investigate the behavior of various methods in a setting where a boundary layer is present in both tangent and normal to the background convection $\bl a$.
This 2D problem, which is shown schematically in Figure \ref{fig:2sch}, is governed in the frequency domain by  
\begin{equation}
\begin{split}
i\omega \phi + a\phi_{,x} & = \kappa (\phi_{,xx} + \phi_{,yy}),\\ 
\phi(x,0) &= \phi(L,y) = 0, \\
\phi(x,L) &= \phi(0,y) = 1. \\
\end{split}
\label{2d-prob}
\end{equation}
Since Dirichlet boundary conditions are imposed on all four boundaries and flow is from left to right, this case produces a flow-facing boundary layer on the right and two flow-tangent boundary layers on top and bottom. 
The boundary layer on the right is similar in nature to what we considered in the 1D problem in Section \ref{sec:1dtest}. 
The other two, however, test the behavior of various methods when there is a sharp gradient in $\phi$ in the crosswind direction, which is a situation that was not tested by that earlier 1D case. 

The selection of this 2D case was motivated by the existence of a closed-form solution that permits us to evaluate the accuracy of all methods. 
As detailed in Appendix \ref{sec:2dsol}, the closed-form solution is obtained in the form of a series using the method of separation of variables and is
\begin{equation}
\begin{split}
\phi(x,y) &= \frac{\sinh\left(\sqrt i W \frac{y}{L}\right)}{\sinh(\sqrt i W)} + \sum_{n=1}^\infty \left[ A_n \exp\left(r_{-n}\frac{x}{L}\right) + B_n \exp\left(r_{+n}\frac{x}{L}\right) \right]\sin\left(\frac{n\pi y}{L}\right), \\
A_n &= \frac{(a_n+b_n)\exp(r_{+n}) - b_n}{\exp(r_{+n}) - \exp(r_{-n})},\;\; B_n = \frac{b_n -(a_n+b_n)\exp(r_{-n})}{\exp(r_{+n}) - \exp(r_{-n})}, \\ 
a_n &= \frac{2(1-\cos(n\pi))}{n\pi},\;\; b_n = \frac{2n\pi\cos(n\pi)}{iW^2 + (n\pi)^2}, \\
r_{\pm n} &= P \pm \sqrt{P^2+iW^2 + (n\pi)^2}. 
\end{split}
\label{2d-sol}
\end{equation}
where the Peclet $P$ and Womersley $W$ number definitions are identical to that of the 1D problem in Eqs. \eqref{P-def} and \eqref{W-def}, respectively. 
In practice, we perform the summations in Eq. \eqref{2d-sol} for 200 terms, which is sufficiently large for the truncation error to be negligible in comparison to the numerical error associated with the methods tested below. 

\begin{figure}
  \centering
  \include{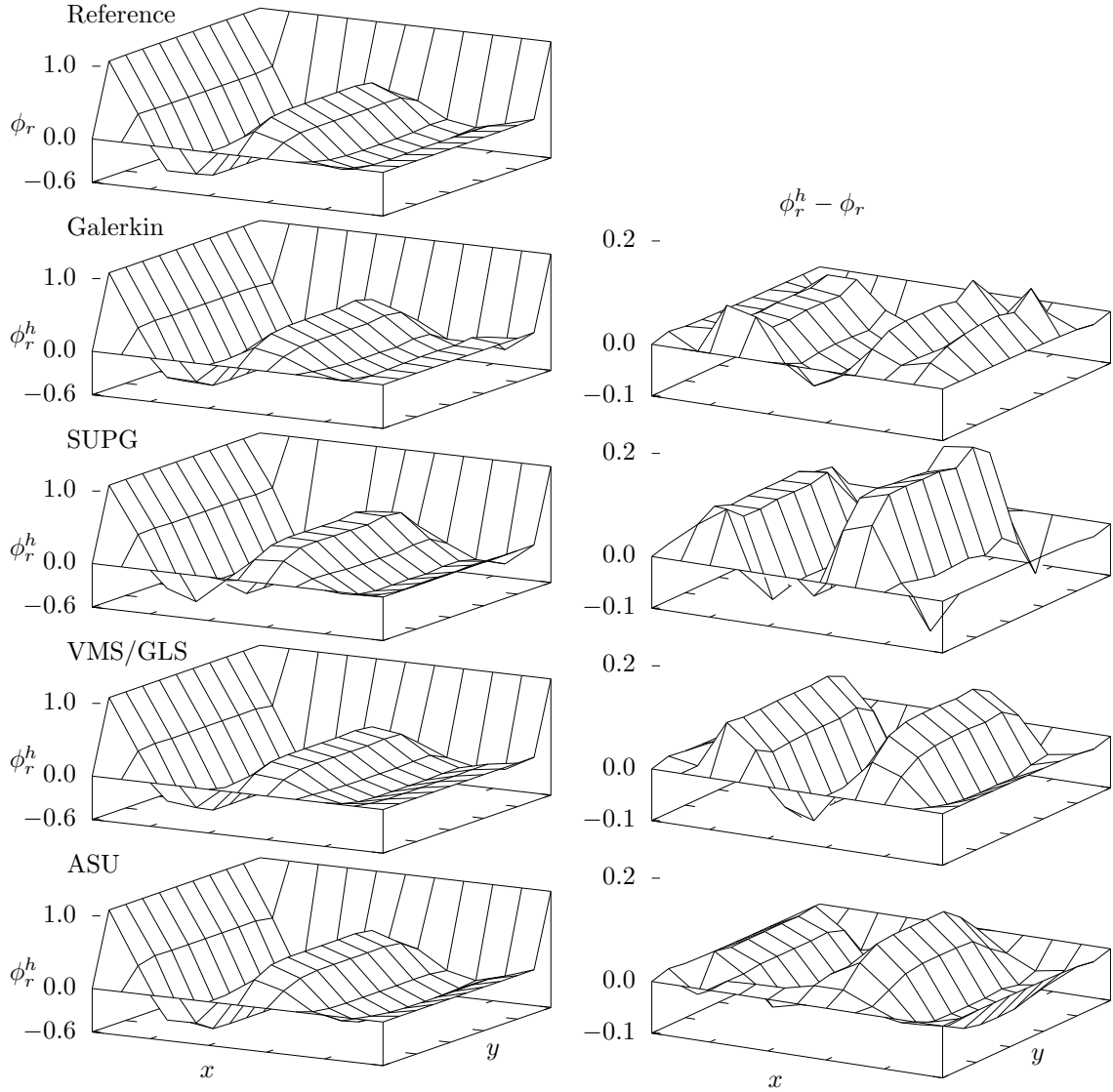}
   \caption{The real part of the reference solution from Eq. \eqref{2d-sol} compared against numerical solutions from Eqs. \eqref{gal-bg}, \eqref{supg-bs}, \eqref{vms-bv}, and \eqref{asu-ba} for the 2D problem shown in Figure \eqref{fig:2sch} at $W=10^{3/2}$ and $P=100/(8\pi)$. The right column is the deviation of numerical solutions from the reference solution.}
  \label{fig:2d} 
\end{figure}

In Figure \ref{fig:2d}, we have compared the reference solution from Eq. \eqref{2d-sol} against the numerical solutions discussed earlier in Section \ref{sec:md}. 
These calculations are performed on a $10\times10$ bilinear grid using $W=10^{3/2}\approx 31.6$ and $P=100/(8\pi)\approx 40$. 
For this specific regime, Galerkin's method performs relatively well as the quickly varying boundary condition (high Womersley number) causes $\phi$ to diffuse before reaching the right boundary, thus avoiding the formation of a boundary layer that could otherwise create non-physical oscillations. 
For the remaining stabilized methods, we observe a trend similar to that of the 1D model problem, with the ASU and SUPG being the best and least accurate methods, respectively, and the VMS/GLS being somewhere in the middle. 

\begin{figure}
  \centering
  \include{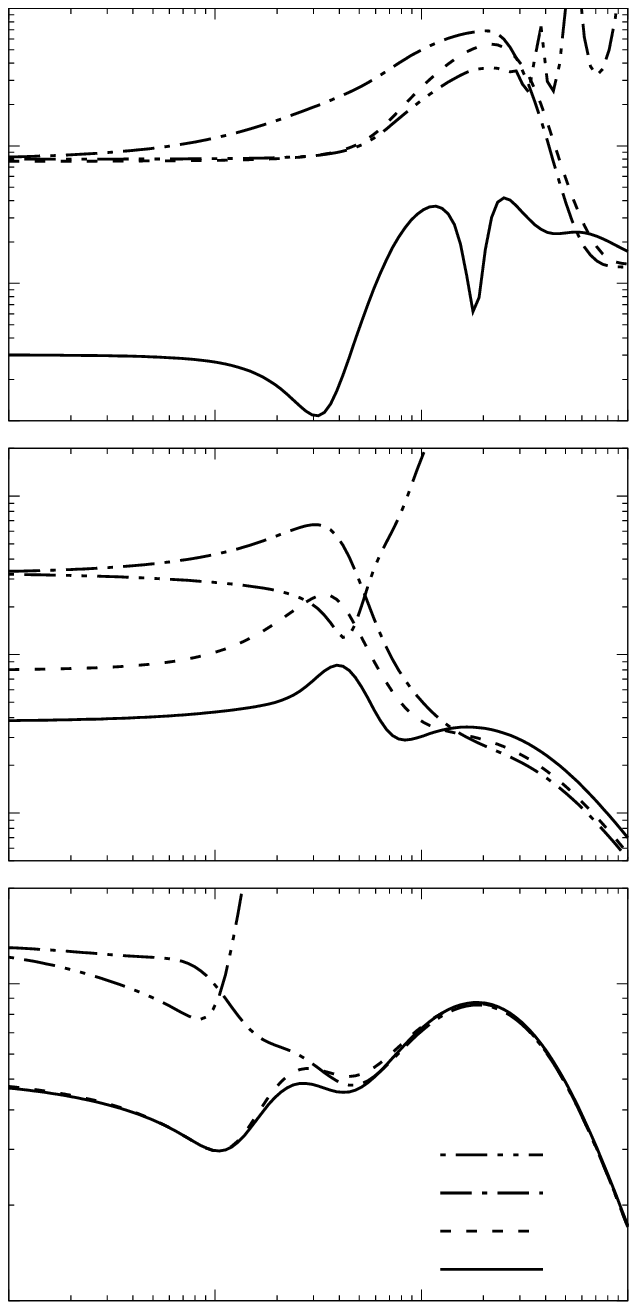}
   \caption{Error in numerical solution of 2D model problem (Figure \ref{fig:2sch}) on a $10\times10$ grid as a function of Peclet number $P$ for Womersley $W=100$ (a), $W=10^{3/2}$ (b) and $W=10$ (c). }
  \label{fig:2err}
\end{figure}

To generalize the above observation, we have repeated these computations using a wider range of conditions. 
The results are presented in terms of $L_2$-norm error in Figure \ref{fig:2err}. 
These results are in agreement with what we discussed earlier, namely the ASU and VMS/GLS methods' similar behavior at small $\beta$ or $W$ and their higher accuracy in comparison to the SUPG at small $\alpha$ or $P$ (Figure \ref{fig:2err}-(c)). 
Also, at large $\alpha$ or $P$, the solutions from all stabilized methods coincide, while that of the Galerkin's method diverges. 
As we have seen with the 1D case, the ASU is the most accurate method for this 2D case as well. 

The observation that the various methods perform similarly considering this 2D case and the earlier 1D case, despite the presence of flow-tangent boundary layers, was somewhat expected. 
Considering the multidimensional discrete form of various methods in Section \ref{sec:md}, the added stabilization terms are either tested by $\tau \bl a.\nabla w^h$ or $\tau \omega w^h$. 
These two sets of terms become important when either $P$ or $W$ are large. 
In the case of a large $P$, where stabilization is required due to the strong convection, $\tau \bl a.\nabla w^h$ is only active in the streamwise direction, thus correctly avoiding the unnecessary introduction of crosswind terms in the discrete form. 
In the case of a large $W$, where stabilization is required due to the fast varying boundary conditions, $\tau \omega w^h$ acts the same in all directions, thus correctly not distinguishing between the flow-facing and flow-tangent boundary layers.  
Note this argument also applies to the ASU as $\hat \omega/\omega$ and $\kappa_{_{\rm ASU}}/\kappa$ are dependent on $\tau\omega$ and $\hat \omega \tau_{\rm diff}$, respectively, which properly act the same in all direction. 

In terms of stability, none of the methods exhibit any issue for the investigated range of conditions. 
For our iterative linear solver, we used the Generalized minimal residual method (GMRES) with a tolerance of $10^{-4}$~\cite{saad1986gmres} that is implemented in our in-house linear solver~\cite{Esmaily2015DS,Esmaily2015BIPN,Esmaily2013PC}. 
The number of GMRES iterations, when averaged over all simulated $P$ above, was larger for Galerkin's method. 
This larger number is caused by the larger number of iterations at higher $P$ where Galerkin's method diverges (Table \ref{table:2d}). 
Among the stabilized methods, the average number of iterations increases with $W$ except for the VMS/GLS method which produces a positive definite matrix. 
This observation is expected as the SUPG and ASU methods can produce a tangent matrix with a wide spectrum of eigenvalues as $\omega$, and correspondingly $W$, increases. 

\begin{table}
\begin{center}
\begin{tabular}{c|ccccc}
 $W$ & Galerkin & SUPG & VMS/GLS & ASU  \\  \hline\hline
    10  & 54 & 20 & 20  & 21 \\
    31.6& 66 & 27 & 22  & 27 \\ 
    100 & 81 & 64 & 21  & 64 \\\hline
    ave.& 67 & 37 & 21  & 37 \\\hline \hline
\end{tabular}
\caption{The average number of GMRES iterations for the 2D cases shown in Figure \ref{fig:2err}.}
\label{table:2d}
\end{center}
\end{table} 

\subsection{A 3D test case}\label{sec:3d}
In the previous case, the edges of the bilinear elements were perfectly aligned with the coordinate directions and that of the flow.
This alignment permits one to directly select a direction-dependent element size $h$ to compute a direction-dependent $\tau$.
Even though such an exercise will produce more accurate results in that idealized setting, it will be hard to generalize those performances to practical cases where there is no clear definition of $h$ for a given direction. 

\begin{figure}
  \centering
  \includegraphics[width=0.45\textwidth]{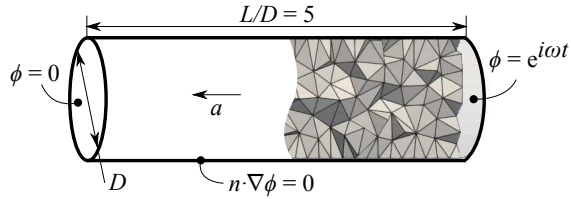}
   \caption{The schematic of the 3D model problem that is defined in a cylindrical domain with a uniform velocity prescribed in the axial direction. A cut of the tetrahedral mesh that is used for all simulations is shown.}
  \label{fig:psch}
\end{figure}

Here, we relied on a scalar definition of $\tau$ that estimated $h$ from the metric tensor (Eq. \eqref{G_def}). 
Nevertheless, to stress-test the methods in a general setting where there is no clear alignment between the physical and parent element coordinate systems, we consider the case shown in Figure \ref{fig:psch}. 

This case involves the unsteady convection-diffusion of a scalar field in a 3D cylindrical domain. 
The imposed boundary conditions and prescribed velocity field are selected to allow for the extraction of a closed-form solution for $\phi$. 
Namely, the velocity is uniform and set $a$ along the cylinder axial direction in the entire domain. 
Dirichlet boundary conditions are imposed at the two ends with one end being zero and the other end being one in the frequency domain. 
A Neumann boundary condition is imposed for the remainder cylindrical shell face, thus effectively reducing this 3D problem to the 1D problem described in Section \ref{sec:1D-problem-statement}.

Even though the exact solution for this 3D problem, which is only a function of position along the axial direction, is identical to the 1D example (namely, Eq. \eqref{1D-exact}), the numerical solution behaves very differently.  
The differentiating factor for the 3D case is that the domain is discretized using tetrahedral elements and the convective velocity does not necessarily align with a parent element coordinate direction. 
As a result, the overall behavior of a given numerical method relies on the multidimensional definition of the stabilization parameter (in particular $\tau$) and whether it correctly captures the effective element size in the streamwise and crosswind directions. 

The specifics of this numerical experiment are as follows. 
The mesh is composed of 2,350 nodes and 8,604 tetrahedral elements, a section of which, is shown in Figure \ref{fig:psch}. 
The length-to-diameter ratio of the cylinder is 5. 
The linear system produced by the various numerical methods is solved using the GMRES iterative method~\cite{saad1986gmres}.
The tolerance on the linear solver is set to $10^{-4}$, which is verified to be sufficiently small so that the reported errors are independent of that tolerance. 
Similar to the 2D test case above, these computations are performed using our in-house finite element solver that has been extensively used for cardiovascular simulations in the past~\cite{lu2022sox17, jia2021efficient, jia2022characterization, Esmaily2012multipleSPS}.

\begin{figure}
  \centering
  \includegraphics[width=1.0\textwidth]{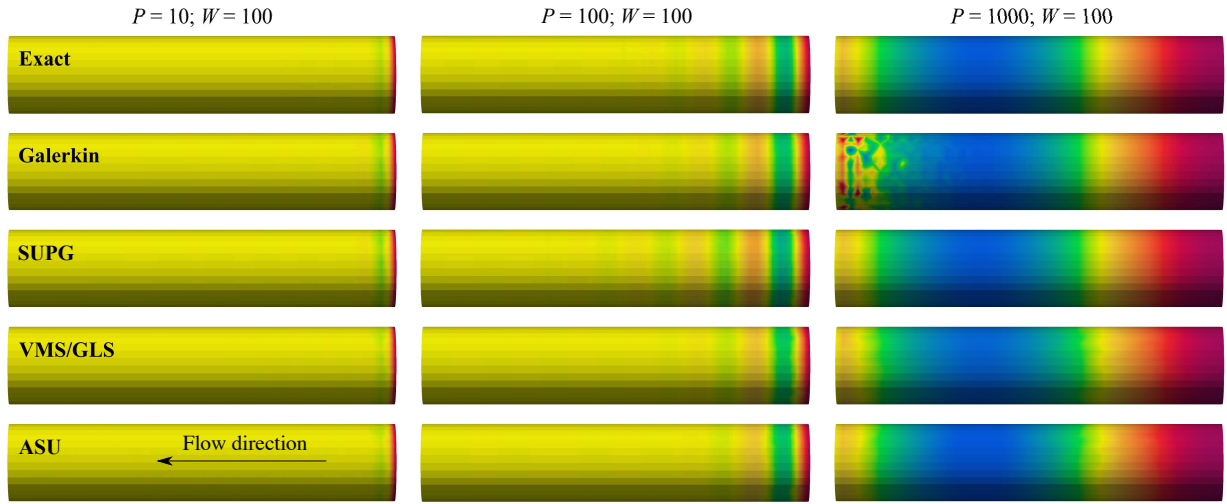}
   \caption{The behavior of the real component of the exact and numerical solutions for the 3D problem shown in Figure \ref{fig:psch} at Womersley number $W=100$ and three values of Peclet number $P=10$ (left column), $P=100$ (middle column), and $P=1000$ (right column).}
  \label{fig:all-pipe}    
\end{figure}

The results of these calculations at three Peclet numbers of $P=10$, 100, and 1000 and a Womersley number of $W=100$ are shown in Figure \ref{fig:all-pipe}. 
Consistent with what we have observed earlier, Galerkin's method generates non-physical oscillations at high $P$. 
All stabilized methods produce reasonably accurate results at small and large $P$, where either the modifications to the baseline Galerkin's method diminish or those modifications are dominated by the SUPG term, respectively. 
At the intermediate $P$, however, the SUPG method exhibits an under-damped behavior by generating oscillations at higher amplitude than the exact solution. 
On the other hand, the ASU and, to a lesser extent, VMS/GLS exhibit over-damped behavior. 
The Galerkin's method produces relatively accurate predictions for this case. 

\begin{figure}
  \centering
  \include{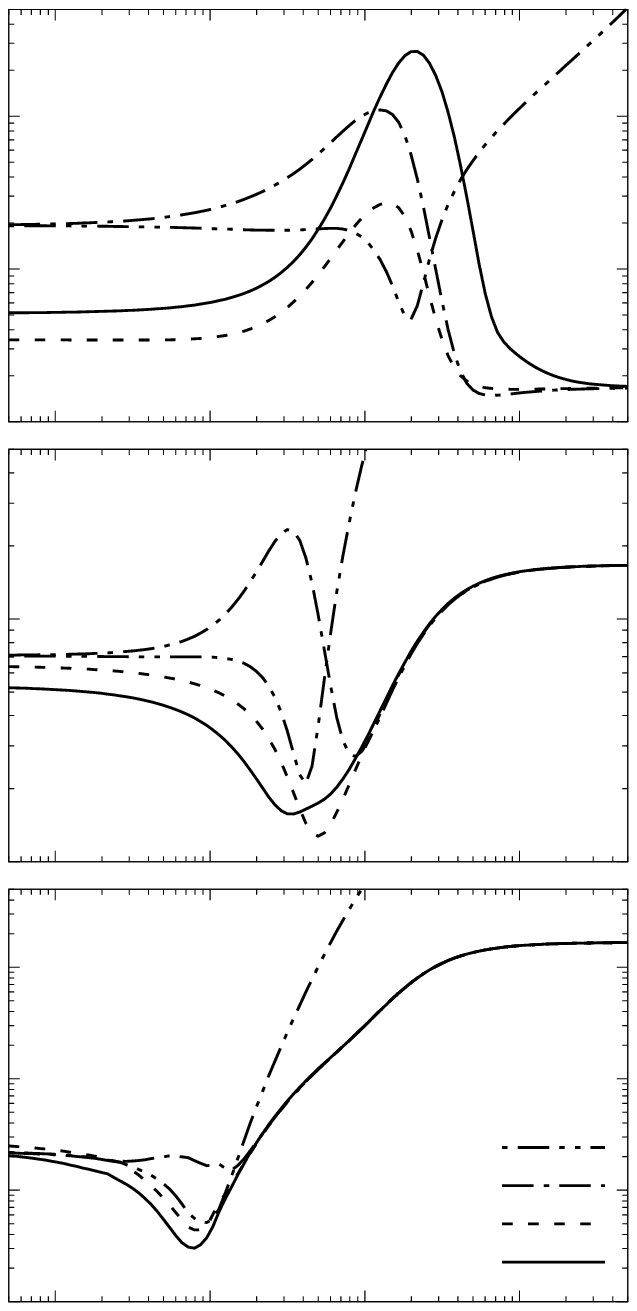}
   \caption{Error in the solution of various methods for the 3D problem shown schematically in Figure \ref{fig:psch} as a function of the Peclet number $P$ at Womersely number of $W=100$ (a; top), $W=31.6$ (b; middle) , and $W=10$ (c; bottom).}
  \label{fig:perr}    
\end{figure}

These computations are repeated for a wider range of $P$ and $W$ and the solution error is reported in Figure \ref{fig:perr} to allow for a more quantitative comparison of the accuracy of all methods. 
The relative behavior of these methods at $W=10$ and $W=31.6$ is similar to what we observed in the 1D and 2D settings. 
The ASU outperforms other techniques, the Galerkin's method diverges at large $P$ and the SUPG is less accurate at smaller $P$. 

For the highly oscillatory $W=100$ condition, a new picture emerges. 
The ASU accuracy is degraded at small $P$ and, to a larger extent, at intermediate $P$. 
As we saw earlier in Figure \ref{fig:all-pipe}, the Galerkin's method produces the most accurate results in those regimes. 
The loss of accuracy of the ASU, in this case, is caused by limiting the $\omega\tau$ angle in the definition of $\hat \omega$ in Eq. \eqref{omegah-lim}.
Owing to the large value of $W=100$, the ASU tangent matrix will become ill-conditioned without that treatment for reasons discussed under Section \ref{sec:conv}. 
With this treatment, the ASU remains stable for all cases considered in Figure \ref{fig:perr}. 

The average number of GMRES iterations for all methods is studied in Table \ref{table:pipe}. 
Since the reported figures are averaged over all simulations with different values of $P$, they reflect outlier simulations that require many linear solver iterations to converge. 
As we saw earlier with the 2D case, the Galerkin's method has a relatively high iteration number at all $W$ owing to the cases with high $P$ where it generates non-physical oscillations. 
The iterative linear solver for the SUPG method converges relatively quickly in all cases, which is somewhat unexpected given that it does not formally produce a positive definite matrix. 
The behavior of the other two stabilized methods is in accordance with what we discussed in Section \ref{sec:conv} regarding the property of their tangent matrices. 
The VMS/GLS method, with a positive definite tangent matrix, exhibits excellent stability. 
While the ASU method remains stable at $W=10$ and $W=31.6$, it struggles at $W=100$, confirming our earlier analysis of this method in Section \ref{sec:conv}. 
Although the linear solver converges for all cases before reaching the set maximum number of iterations of 1,000, it requires a larger number of iterations on average, indicating the relative instability of the ASU method. 

\begin{table}
\begin{center}
\begin{tabular}{c|cccc}
 $W$ & Galerkin & SUPG & VMS/GLS & ASU  \\  \hline\hline
    10  & 160 & 102& 102 & 105 \\
    31.6& 129 & 71 & 68  & 81  \\ 
    100 & 128 & 69 & 44  & 292 \\\hline
    ave.& 139 & 81 & 71  & 159 \\\hline \hline
\end{tabular}
\caption{The average number of GMRES iterations for the 3D cases shown in Figure \ref{fig:perr}.}
\label{table:pipe}
\end{center}
\end{table} 

\subsection{Mesh-convergence study} \label{sec:msh}
In Appendix \ref{app:cons}, we showed that all the methods under consideration are at least second-order accurate for the 1D model problem as $h\to0$. 
Here, the generalization of those results to multiple dimensions is discussed using numerical experiments and existing error estimates from the literature \cite{johnson1984finite,Hughes1986beyond}.

To investigate the behavior of the error versus mesh size $h$, we consider the 3D case from Section \ref{sec:3d}. 
The simulations at two values of Peclet number $P=10$ and $P=1000$ and two values of Womersely number $W=10$ and $W=100$ are performed using various methods while the mesh is being successively refined from 2,350 to over a million elements. 
All remaining parameters are identical to the case in Section \ref{sec:3d} except for the linear solver tolerance that is reduced to $10^-{12}$ to ensure errors at finer grids are independent of that parameter.

\begin{figure}
  \centering
  \include{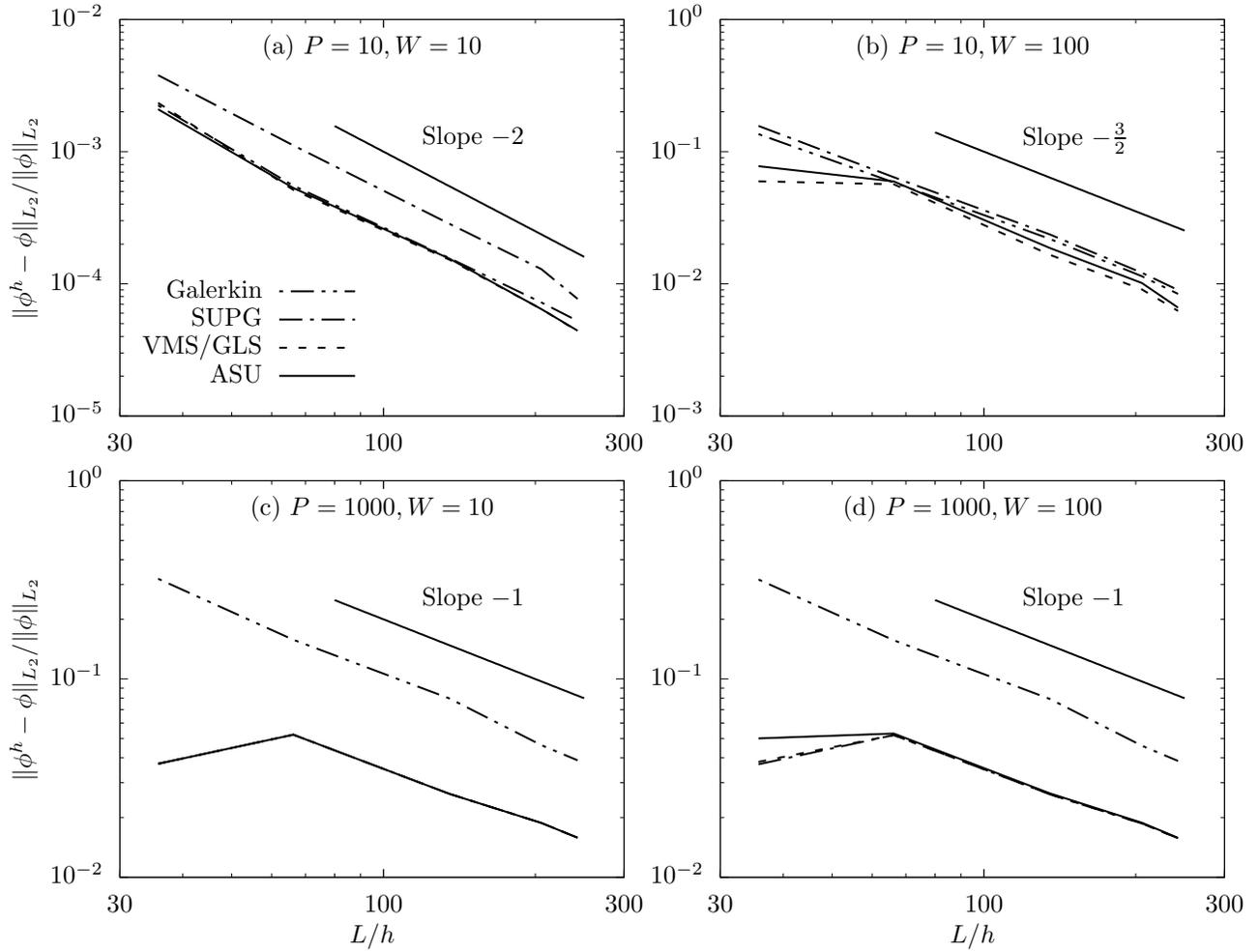}
   \caption{Error in the solution of various methods for the 3D problem shown schematically in Figure \ref{fig:psch} as a function of element size $h$ (normalized by the pipe length $L$) at the  Womersely number of $W=100$ for the Peclet number $P=1000$ (a; left) and $P=10$ (b; right).}
  \label{fig:pmsh}    
\end{figure}

The results of these calculations is shown in Figure \ref{fig:pmsh}. Considering the value of $P$ and $W$ and the simulated range of $L/h$, Figure \ref{fig:pmsh}-(a) falls in the diffusive limit where $\alpha\ll1$, Figure \ref{fig:pmsh}-(b) falls in the oscillatory regimes where $\beta$ is of order 1, and lastly, Figure \ref{fig:pmsh}-(c) and (d) fall in the convective limit where $\alpha> 1$.  
A few observations can be made based on these results:
\begin{enumerate}
\item Regardless of the regime under consideration, all methods show the same order of convergence. 
\item In the diffusive limit (Figure \ref{fig:pmsh}-(a)), all methods are second order. This observation is compatible with what we found in Appendix \ref{app:cons} and the existing estimates in the literature for the diffusive limit showing $\|\phi^h-\phi\|_{L_2} = O(h^{k+1})$ with $k=1$ for linear elements \cite{Hughes1986beyond}.
\item In the convective limit, the order of convergence drops by one. This sub-optimal order of convergence is compatible with the existing error estimates in the literature. For the baseline Galerkin's method $\|\phi^h-\phi\|_{L_2} = h^k\|\phi\|_{H^{k+1}}$ has been provided, whereas for streamline diffusion methods (e.g., SUPG) better error scaling of $\|\phi^h-\phi\|_{L_2} = h^{k+1/2}\|\phi\|_{H^{k+1}}$ has been reported~\cite{johnson1984finite, Hughes1986beyond} for stationary problems. The observed scaling in Figure \ref{fig:pmsh}-(c) and (d), which is obtained for non-stationary problems, is closer to the former at $O(h)$ than the latter estimate at $O(h^{3/2})$. The latter estimate, however, is observed at the highly oscillatory regime shown in Figure \ref{fig:pmsh}-(b). 
\end{enumerate}

\subsection{Patient-specific case and comparison against time formulation}\label{sec:patient}

This patient-specific case is adopted to evaluate the performance of time-spectral methods in an engineering-relevant setting and compare them against the existing methods that are formulated in time.

\begin{figure}
    \centering
    \includegraphics[width=0.5\textwidth]{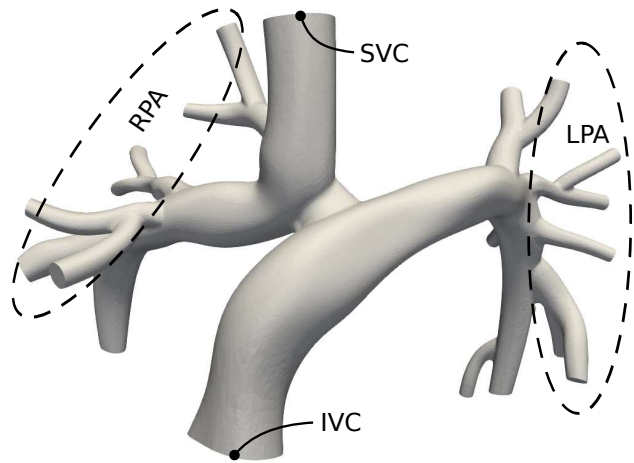}
    \caption{The patient-specific Fontan geometry. }
    \label{fig:fontan}
\end{figure}
The simulated anatomy belongs to a 3-year-old male patient who has undergone Fontan operation~\cite{wilson2013vascular, marsden2010new}.
As shown in Figure \ref{fig:fontan}, the reconstructed geometry contains the superior and inferior vena cava (SVC and IVC, respectively) that redirect de-oxygenated blood to the right and left pulmonary arteries (RPA and LPA, respectively) branches through a cross-shaped connection. 
Clinically, it is important for the hepatic flow from the IVC to be equally distributed to the RPA and LPA to prevent pulmonary arteriovenous malformations \cite{duncan2003pulmonary, yang2013optimization}. 
Computationally, this flow distribution can be obtained by modeling the convection-diffusion of a neutral tracer $\phi$ that is released through the IVC, and subsequently, measuring the flux of $\phi$ through the LPA and RPA. 

Given the clinical relevance of this problem, we simulated the release of a neutral tracer through the IVC using various time-spectral methods discussed in this article. 
Additionally, we compare the results against the standard stabilized finite element methods for the convection-diffusion equation that -- in the absence of a volumetric source and nonzero Neumann boundary -- is formulated as 
\begin{equation}
    (w^h, \hat \phi_{,t}^h) + (w^h,\bl a\cdot  \nabla \hat \phi^h) + (\nabla w^h,\kappa \nabla \hat \phi^h) +\sum_e \Big(\tau \bl a\cdot \nabla w^h, \hat r(\hat \phi^h) \Big)_{\Omega_e} = 0,
\label{supg-time}
\end{equation}
where $\hat \phi^h(\bl x,t)$ is the concentration of tracer at location $\bl x$ and time $t$, and
\begin{equation}
 \hat r(\hat \phi^h) = \hat \phi_{,t} + \bl a\cdot  \nabla \hat \phi^h - \nabla\cdot(\kappa \nabla \hat \phi^h),
\label{res-time}
\end{equation}
is the residual of the temporal form of the convection-diffusion equation. 
Since $\tau$ in Eq. \eqref{supg-time} is defined based on Eq. \eqref{tau_supg-md}, this temporal formulation in the steady-state limit simplifies to the stabilized methods discussed in Section \ref{sec:md}. 
For unsteady cases, this formulation corresponds to the SUPG method (Eq. \eqref{supg-bs}) if one neglects the temporal discretization error.
Thus, by using this formulation, we can perform an apple-to-apple comparison between the temporal and spectral formulations in a realistic setting.

Equation \eqref{supg-time} is discretized in time using the second-order generalized-$\alpha$ time integration scheme \cite{Jansen2000305} with $\rho_\infty=0$ , which effectively simplifies it to the Crank-Nicolson time integration scheme. 
The baseline steady blood flow $\bl a(\bl x)$ is simulated using our exiting stabilized fluid solver \cite{Esmaily2015BIPN}. 
For the flow simulations, we use use physiologic flow rate of 39.9 mL/s that is divided equally and imposed through a parabolic profile at the IVC and SVC. 
Zero Dirichlet and Neumann boundary conditions are imposed on the wall and all outlets, respectively. 
The geometry is discretized using approximately 1.7 million tetrahedral elements. 
Blood viscosity and density are set to physiological values of 0.04 g.cm/s and 1.06 g/cm$^3$, respectively.

To perform this comparison in an unsteady regime, the concentration of the tracer at the IVC is taken to vary periodically during one cardiac cycle. 
Using a heart rate of 60 beats per minute, $\hat \phi = \cos(2\pi t)$ is imposed uniformly at the IVC. Also, a zero Dirichlet boundary condition is imposed at the SVC. Zero Neumann boundary conditions are imposed for the remaining boundaries.
The molecular diffusivity of the tracer is set to $\kappa = 10^{-2}$ cm$^2$/s. 
For the temporal simulation, we use 500 time steps per cardiac cycle to ensure temporal discretization error is sufficiently small. 
Three cardiac cycles are simulated to wash out the zero initial condition, requiring a total of 1,500 time steps. 
The same linear solver (GMRES) with the same tolerance of $10^{-6}$ is used for both temporal and spectral methods. 
The computations are performed using 32 compute cores on a single node with dual Xeon E5-2698 v4 processors.

\begin{figure}
    \centering
    \includegraphics[width=1\textwidth]{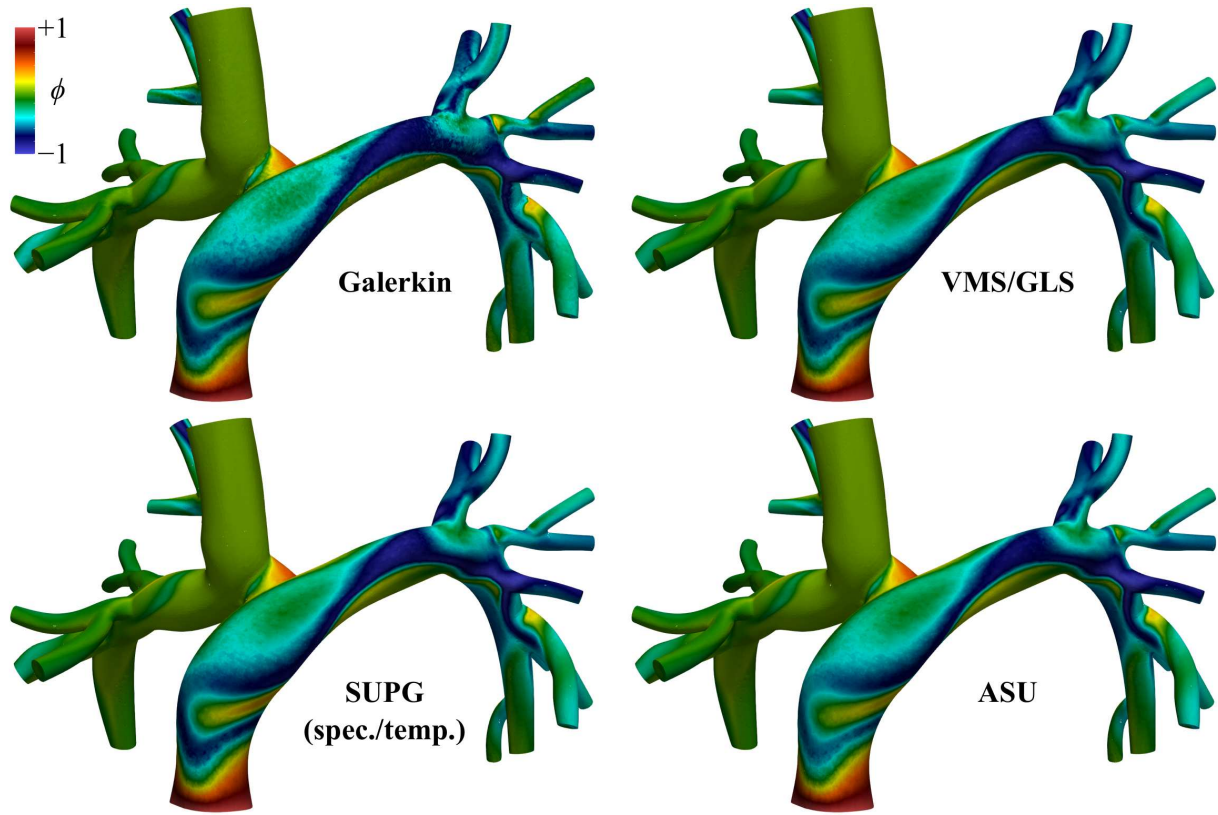}
    \caption{The simulated concentration of a tracer, which is released periodically through the IVC in the domain, using various methods.}
    \label{fig:fontan-sol}
\end{figure}
The result of these calculations is shown in Figure \ref{fig:fontan-sol} for various methods. 
These results correspond to $\phi^h_r$, which is the same as $\hat \phi^h$ at the beginning of the cardiac cycle. 
Since the temporal and spectral SUPG results were almost identical (relative difference of $\approx 10^{-4}$), only the spectral SUPG result is shown in this figure. 
The general behavior of the solution is the same as the earlier cases with the Galerkin's method producing some spurious oscillations that are eliminated by the stabilized methods.

The computational performance of various methods is studied in terms of the number of linear solver iterations and solution turnover time in Table \ref{table:fontan}. 
A few observations can be made based on these results:
\begin{enumerate}
    \item All the stabilized time-spectral methods perform similarly, producing a solution in less than 10 seconds for this mesh that has over 1 million elements. 
    \item The linear solver struggles to converge for the time-spectral Galerkin's method that produces a highly spurious solution (Figure \ref{fig:fontan-sol}) in this convection-dominated regime.
    \item The conventional SUPG method that is formulated in time is approximately 40 times more expensive than its spectral counterpart. This gap in performance is explained by the fact that the temporal formulation requires over 40,000 linear iterations in total for 1,500 time steps, whereas the spectral formulation requires less than 1,000. In general, the performance gap between the spectral and temporal formulations depends on the number of simulated modes, the number of time steps per cardiac cycle, and the number of simulated cycles to wash out the initial conditions. 
\end{enumerate}
\begin{table}
\begin{center}
\begin{tabular}{c|cc}
 Method   & GMRES iterations & Simulation time (sec.)  \\  \hline\hline
 Galerkin & 6,282 & 54 \\
 SUPG     & 853  & 7.8 \\ 
 VMS/GLS  & 831  & 7.2 \\
 ASU      & 861  & 7.4 \\\hline 
 Time-SUPG& 41,986$^*$ & 290 \\\hline\hline
\end{tabular}
\caption{The number of iterative solver (GMRES) iterations and simulation turnover time for the simulations shown in Figure \ref{fig:fontan-sol}. $^*${\small This is the cumulative number of iterations from all 1,500 time steps.} }
\label{table:fontan}
\end{center}
\end{table} 


\section{Conclusions} \label{sec:conclusions}
We introduced and compared five methods for the solution of the time-spectral convection-diffusion equation. 
These methods are evaluated based on accuracy, stability, and convergence properties. 
The baseline Galerkin's method (Eq. \eqref{gal}) for the time-spectral equation behaves very similarly to its conventional time formulation counterpart, producing non-physical oscillations in strongly convective regimes. 
Including the streamline upwind Petrov/Galerkin (SUPG) stabilization term in the formulation removes those non-physical oscillations. 
Nevertheless, the resulting method (Eq. \eqref{supg}) tends to produce a solution that overshoots physical oscillations at high Womersley numbers. 
We also explored the variational multiscale (VMS) and Galerkin/least-square (GLS) methods, which are identical for the problem under consideration (Eq. \eqref{vms}).
The VMS/GLS method demonstrated reasonable accuracy across all cases studied.
Moreover, it produces a formally positive-definite tangent matrix, leading to excellent stability as verified by our numerical experiments.
The last method introduced was the augmented SUPG (ASU), which we tailor-designed to produce a nodally exact solution for the time-spectral convection-diffusion equation in 1D. 
This method (Eq. \eqref{asu}) achieves the highest accuracy across all tested cases, except for the 3D case at the highest Womersley number, where a modification (Eq. \eqref{omegah-lim}) is necessary for stability.
In terms of mesh convergence, all methods behaved the same, showing a first and second-order convergence with the element size in the convective and diffusive limits, respectively.
Considering all factors, the VMS method is an attractive option for achieving a balance between accuracy, stability, and implementation simplicity. 
The ASU, on the other hand, is an optimal choice for regimes with low to moderate element Womersley number ($\beta\le O(1)$) if accuracy takes precedence. 

\section{Acknowledgement}
For the patient-specific case in Section~\ref{sec:patient}, the geometry used herein was provided in whole or in part with Federal funds from the National Library of Medicine under Grant No. R01LM013120, and the National Heart, Lung, and Blood Institute, National Institutes of Health, Department of Health and Human Services, under Contract No. HHSN268201100035C.
\bibliographystyle{unsrt}
\bibliography{main}


\newpage
\appendix
\section{Galerkin's solution to the 1D problem} \label{appendix1}

Galerkin's solution $\phi^h(x)$ is expressed in terms of piecewise linear shape function, which for a given node $A$ is expressed as 
\begin{equation}
    N_A(x) = \left\{ 
    \begin{matrix}
        \frac{x-x_{A-1}}{x_A - x_{A-1}},  & x_{A-1}<x<x_A \\
        \frac{x_{A+1}-x}{x_{A+1} - x_A}, & x_A<x<x_{A+1} \\
        0, & {\rm otherwise}.  
    \end{matrix} 
    \right.
    \label{A1}
\end{equation}
Using these shape functions, the solution and test functions are expressed as
\begin{equation}
    \phi^h(x) =  \sum_{A=0}^{N} U_A N_A(x) , 
    \label{A2}
\end{equation}
and 
\begin{equation}
    w^h(x) = \sum_{A=1}^{N-1} c_A N_A(x), 
    \label{A3}
\end{equation}
respectively, where $U_A$ is the solution at node $A$ and $c_A$ is an arbitrary constant associated with node $A$. 
Note that the end nodes $A=0$ and $N$ are excluded from the summations in Eq. \eqref{A3} as Dirichlet boundary conditions are imposed on the two ends of the domain. 
Also, since $\phi^h(0)=0$ and $\phi^h(L)=1$, we require $U_0=0$ and $U_N=1$ in what follows. 

Using Eqs. \eqref{A2} and \eqref{A3} to plug in for $\phi^h$ and $w^h$ in Eq. \eqref{1D-weak} produces
\begin{equation}
    \sum_{A=1}^{N-1} c_A \left\{ \sum_{B=0}^{N} \left[ \bigg(N_A,i\omega N_B\bigg)  + \left(N_A,a\frac{\partial N_B}{\partial x}\right) +\left(\frac{\partial N_A}{\partial x},\kappa \frac{\partial N_B}{\partial x}\right) \right]U_B \right\}  = 0,
\label{A4}
\end{equation}

Since Eq. \eqref{A4} must hold for any $c_A$, we can conclude 
\begin{equation}
\sum_{B=0}^{N} \left[\bigg(N_A,i\omega N_B\bigg)  + \left(N_A,a\frac{\partial N_B}{\partial x}\right) +\left(\frac{\partial N_A}{\partial x},\kappa \frac{\partial N_B}{\partial x}\right)\right]U_B  = 0, \;\;\; A = 1, \cdots N-1. 
\label{A5}
\end{equation}

The three inner products in Eq. \eqref{A5} can be calculated explicitly using the fact that we are interested in a uniform grid of spacing $h$ and $x_{A+1}-x_A= x_A-x_{A-1} = h$ in Eq. \eqref{A1}.   
From Eq. \eqref{A1} it is straightforward to show 
\begin{equation}
    (N_A,N_B) = \left\{ 
    \begin{matrix}
        \frac{1}{6}h,  & B=A\pm 1 \\
        \frac{2}{3}h, & B=A \\
        0, & {\rm otherwise}.  
    \end{matrix} 
    \right.
    \label{A6}
\end{equation}
similarly
\begin{equation}
    \left(N_A,\frac{\partial N_B}{\partial x}\right) = \left\{ 
    \begin{matrix}
        \pm \frac{1}{2},  & B=A\pm 1 \\
        0, & {\rm otherwise},  
    \end{matrix} 
    \right.
    \label{A7}
\end{equation}
and
\begin{equation}
    \left(\frac{\partial N_A}{\partial x}, \frac{\partial N_B}{\partial x}\right) = \left\{ 
    \begin{matrix}
        -\frac{1}{h},  & B=A\pm 1 \\
        \frac{2}{h}, & B=A \\
        0, & {\rm otherwise}.  
    \end{matrix} 
    \right.
    \label{A8}
\end{equation}

Plugging in for the inner products that appear in Eq. \eqref{A5} from Eqs. \eqref{A6}, \eqref{A7}, and \eqref{A8} yields
\begin{equation}
    \left(\frac{i\omega h}{6} + \frac{a}{2} -\frac{k}{h}\right)U_{A+1} + \left(\frac{2i\omega h}{3} + \frac{2\kappa}{h}\right)U_A + \left( \frac{i\omega h}{6} - \frac{a}{2} -\frac{k}{h} \right)U_{A-1} = 0.
    \label{A9}
\end{equation}

Multiplying Eq. \eqref{A9} by $h/\kappa$ will non-dimensionalize the coefficients, allowing us to express them in terms of $\alpha$ and $\beta$ from Eqs. \eqref{alpha} and \eqref{beta}, respectively.  
The result is
\begin{equation}
    (i\beta + \alpha - 1)U_{A+1} + (4i\beta + 2)U_A + ( i\beta - \alpha - 1)U_{A-1} = 0.
    \label{A10}
\end{equation}

Provided that the exact solution has an exponential form (Eq. \eqref{1D-exact}), it is reasonable to assume that the numerical solution also takes an exponential form. 
The exponent, however, can differ from that of the exact solution. 
That difference can be accommodated by selecting an exponent base that is a free parameter, which by following the tradition we call $\rho$. 
That permits us to take 
\begin{equation}
     U_A = c\rho^A, 
     \label{A11}
\end{equation}
in which the arbitrary constant $c$ is included to allow for matching the solution to the given boundary conditions. 

For the solution guess from Eq. \eqref{A11} to satisfy Eq. \eqref{A10}, $\rho$ must be the roots of a quadratic polynomial that is 
\begin{equation}
    (i\beta + \alpha - 1)\rho^2 + (4i\beta + 2)\rho + i\beta - \alpha - 1 = 0.
    \label{A12}
\end{equation}
The two roots of Eq. \eqref{A12} are those that were shown in Eq. \eqref{rho12}. 

The general solution can be expressed in terms of the two discrete solutions associated with the two roots $\rho_{1,2}$, namely 
\begin{equation}
    U_A = c_1\rho_1^A + c_2\rho_2^A,
    \label{A13}
\end{equation}
where $c_1$ and $c_2$ are arbitrary constants selected so that the two boundary conditions are satisfied. 
Since $U_0  = 0$ and $U_N = 1$, it is straightforward to show 
\begin{equation}
    c_{1,2} = \pm\frac{1}{\rho_1^N-\rho_2^N}.
    \label{A14}
\end{equation}
Combining Eqs. \eqref{A13} and \eqref{A14} and noting $U_A = \phi^h(x_A)$ produces the result shown in Eq. \eqref{1d-gal}, thus completing this proof.

\section{\cola{Reconstruction of the diffusive flux}}
\label{app:CSUPG}

\cola{The SUPG and GLS/VMS stabilized methods that were discussed in Section \ref{sec:1d} rely on the residual $r_r(\phi^h)$ and $r_i(\phi^h)$ in their construct. When evaluating these residuals in the element interior, the diffusive term $\kappa \phi^h_{,xx}$ drops owing to the fact that $\phi^h$ was discretized using linear interpolation functions. Our goal here is to compare those methods against methods in which that diffusive term in the residual is estimated. 
Since that requires reconstruction of diffusive terms, we refer to these alternative methods as RD-SUPG and RD-VMS in what follows. } 

\cola{To estimate the second derivative of $\phi^h$ at the Gauss point, we follow the method described in \cite{jansen1999better}. 
The overall idea is first to compute and project $\phi^h_{,x}$ to the mesh nodes to obtain a new variable $\psi^h$. Then use $\psi^h$ to estimate $\phi^h_{,xx} = \psi^h_{,x}$ at the Gauss points. 
In practice, this process involves solving}
\begin{equation}
    \cola{(w^h,\psi^h) = (w^h,\phi^h_{,x}),}
    \label{AA1}
\end{equation}
\cola{as the first step. 
The resulting variable $\psi^h$, which is available at the mesh nodes, is continuous across elements and can be differentiated to compute $\psi^h_{,x}$ at the Gauss points to estimate $\phi^h_{,xx}$.
It is this reconstructed $\phi^h_{,xx}$ that is used in the residuals of Eqs. \eqref{1D-SUPG} and \eqref{1D-vms} to obtain RD-SUPG and RD-VMS methods, respectively. }

\begin{figure}[H]
  \centering
  \include{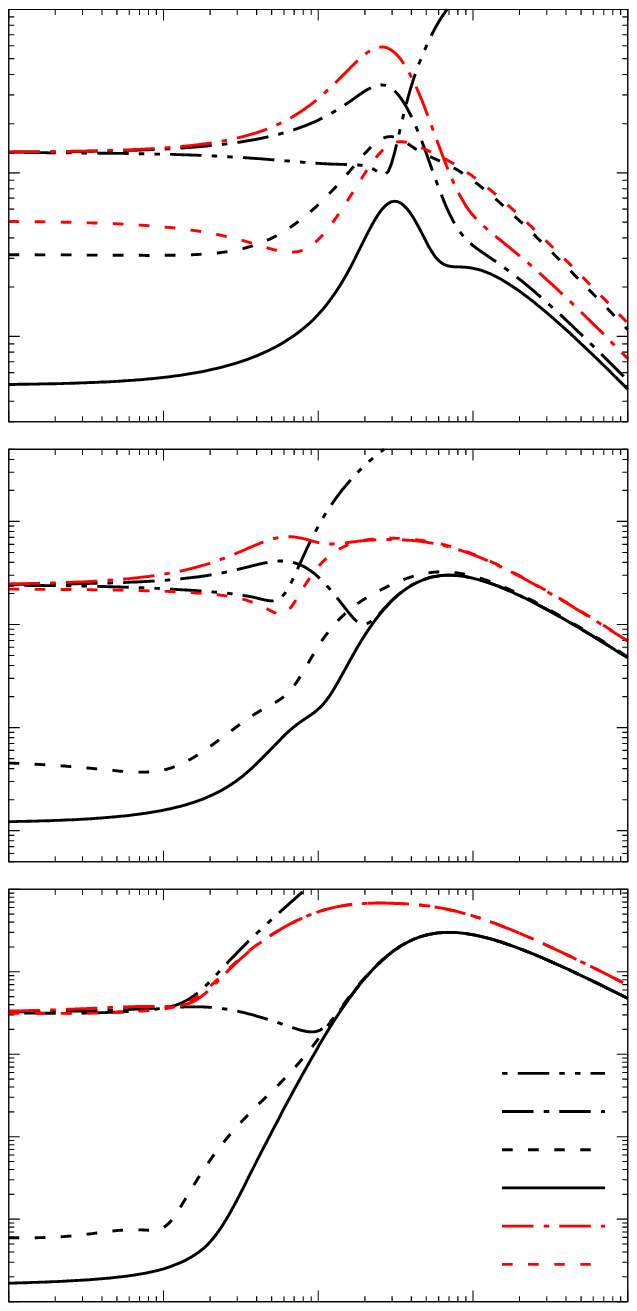}
   \caption{\cola{Same as Figure \ref{fig:err} but also including the error for methods with reconstructed diffusion term in the residual, namely RD-SUPG (red dash-dotted) and RD-VMS (red dotted).} }
  \label{fig:errC}    
\end{figure}

\cola{To compare RD-SUPG and RD-VMS against the baseline SUPG and VMS methods, the test case described under Section \ref{sec:1dtest} is repeated to reproduce Figure \ref{fig:err} for all methods. 
The result of that process is shown in Figure \ref{fig:errC}. }

\cola{Remarks:}
\begin{enumerate}
    \item \cola{The baseline methods that are considered throughout this study perform better than those with reconstructed diffusive flux. This behavior was expected given that the derivation of the standard SUPG method in Section \ref{sec:1dsupg} was based on a closed-form solution that excluded the reconstructed diffusion term. In other words, the SUPG method is nodally exact by design for steady flows when the diffusive term is dropped from the residual. The error is no longer expected to be optimal when that baseline method is modified by incorporating the diffusive flux to obtain RD-SUPG.}
    \item \cola{Solving Eq. \eqref{AA1} as shown requires solving a linear system with a mass matrix on the left-hand side. To avoid this relatively expensive operation, one may adopt a lumped mass matrix instead \cite{Hughes2000finite}. Doing so will generate results that are slightly less accurate than those shown in Figure \ref{fig:errC}. Hence, we opted to show the version that is more accurate but also more expensive.}
    \item \cola{The RD-GLS method was not discussed because its implementation is not straightforward for it requires the reconstruction of the second derivative in the weight space. Incorporating the diffusive flux in the solution space only will reduce that method to RD-VMS, which is already investigated here. }
\end{enumerate}

\section{The derivation of Eqs. \eqref{alphah1} and \eqref{betah1}} \label{app:pres}
For a nodally exact solution, we must have 
\begin{equation}
    \hat \rho_{1,2}^A = \exp(r_{1,2}\frac{x_A}{L}), 
\label{B1}
\end{equation}
where $r_{1,2}$ are computed from the exact solution in Eq. \eqref{r12} and 
\begin{equation}
    \hat \rho_{1,2} = \frac{1 + 2i\hat \beta \pm \sqrt{\hat \alpha^2 - 3\hat \beta^2 + 6i\hat \beta}}{1 - \hat \alpha - i\hat \beta}, 
\label{B2}
\end{equation}
which is the same as Eq. \eqref{rho12} but with the adjusted $\hat \alpha$ and $\hat \beta$ of the proposed augmented SUPG method from Eqs. \eqref{alphah} and \eqref{betah}. 

Since $x_A/L = A h/L$, the right hand side of Eq. \eqref{B1} can be written as 
\begin{equation}
    \exp(r_{1,2}\frac{x_A}{L}) = \left(\exp(r_{1,2}\frac{h}{L})\right)^A. 
\label{B3}
\end{equation}
Thus, from Eqs. \eqref{B1} and \eqref{B3}, we have
\begin{equation}
    \hat \rho_{1,2} = \exp(r_{1,2}\frac{h}{L}).
\label{B4}
\end{equation}
Substituting for $r_{1,2}$ using Eq. \eqref{r12} and simplifying the result by using the fact that $Ph/L = \alpha$ and $(Wh/L)^2 = 6\beta$ (see Eqs. \eqref{P-def}, \eqref{W-def}, \eqref{alpha}, and \eqref{beta}) produces
\begin{equation}
    \hat \rho_{1,2} = \exp(\alpha \pm \sqrt{\alpha^2 + 6i\beta}),
\label{B5}
\end{equation}
or alternatively
\begin{equation}
\begin{split}
    \hat \rho_1 + \hat \rho_2 &= 2\exp(\alpha) \cosh(\gamma),\\
    \hat \rho_1 - \hat \rho_2 &= 2\exp(\alpha) \sinh(\gamma),
\end{split}    
\label{B6}
\end{equation}
where $\gamma$ is defined in Eq. \eqref{gamma}. 
Combining Eqs. \eqref{B2} and \eqref{B6}, yields
\begin{equation}
    \frac{1 + 2i\hat \beta}{1 - \hat \alpha - i\hat \beta} = \exp(\alpha) \cosh(\gamma),
\label{B7}
\end{equation}
and
\begin{equation}
    \frac{\sqrt{\hat \alpha^2 - 3\hat \beta^2 + 6i\hat \beta}}{1 - \hat \alpha - i\hat \beta} = \exp(\alpha) \sinh(\gamma).
\label{B8}
\end{equation}
From Eq. \eqref{B7}, it is straightforward to show
\begin{equation}
   \hat \alpha = 1 - i\hat \beta  - \frac{1 + 2i\hat \beta}{\exp(\alpha) \cosh(\gamma)}. 
\label{B9}
\end{equation}
Substituting for $\hat \alpha$ from Eq. \eqref{B9} into Eq. \eqref{B8} and simplifying the results produces the following quadratic relationship for $i\hat \beta$
\begin{equation}
   \left[4\cosh(\alpha) + 2\cosh(\gamma)\right](i\hat \beta)^2 + \left[4\cosh(\alpha) -  \cosh(\gamma)\right] i\hat \beta + \cosh(\alpha) -\cosh(\gamma)=0. 
\label{B10}
\end{equation}
Equation \eqref{B10} has two roots. 
One root is $i\hat \beta = -1/2$, which is not admissible as it produces zero divided by zero in Eq. \eqref{B7}. 
It is the second root that in combination with Eq. \eqref{B9} produces Eqs. \eqref{alphah1} and \eqref{betah1}, thus completing the proof.

\section{Derivation of Eq. \eqref{omegah}} \label{omegah-app}
To arrive at Eq. \eqref{omegah} from \eqref{omegahe} for regimes in which $\beta \lessapprox 1$, we investigate the asymptotic behavior of $\hat \omega$ at limits of $\alpha \ll 1$ and $\alpha \gg 1$. 

Let us first consider $\alpha \gg 1$. 
Since in this limit $\beta/\alpha^2 \ll 1$, Eq. \eqref{gamma} can be expressed as 
\begin{equation}
\gamma = \alpha + 3i\left(\frac{\beta}{\alpha}\right) + \left(\frac{9\beta^2}{2\alpha^3}\right) + \mathcal O(\alpha^{-5}). 
    \label{C1}
\end{equation}
Since $\exp(\alpha) \gg \exp(-\alpha)$, we have
\begin{equation}
\cosh(\gamma) = \frac{1}{2}\exp(\alpha)\left[1 + \left(\frac{3i\beta}{\alpha}\right) - \left(\frac{9\beta^2}{2\alpha^2}\right) -\left(\frac{9i\beta^3}{2\alpha^3}\right) + \left(\frac{9\beta^2}{2\alpha^3}\right)+ \mathcal O(\alpha^{-4}) \right]. 
    \label{C2}
\end{equation}
Using this result along with the fact that $\cosh(\alpha)\approx \sinh(\alpha) \approx \exp(\alpha)/2$ permits us to simplify Eq. \eqref{omegahe} and obtain
\begin{equation}
\hat \omega = \left(1 + \frac{3i\beta}{2\alpha} - \frac{3\beta^2}{2\alpha^2}\right) \omega + \mathcal O(\alpha^{-2}), \;\;\; {\rm for} \;\;\; \alpha\gg 1.
    \label{C3}
\end{equation}
In deriving Eq. \eqref{C3}, we dropped $-3i\beta/(2\alpha^2)$ in comparison to the $3i\beta/(2\alpha)$. 
The term $-3\beta^2/(2\alpha^2)$, however, was not dropped at it was the leading order in excess of 1 and it determines the in-phase contribution of the stabilization term at large $\alpha$.  

Next, let us consider $\alpha \ll 1$. 
By selecting $\alpha$ such that $\alpha \ll \beta$, we can write 
\begin{equation}
\gamma = \sqrt{6i\beta}\left(1 + \frac{\alpha^2}{12i\beta} + \mathcal O(\alpha^{4})\right). 
    \label{C4}
\end{equation}
Performing Taylor series expansion of $\cosh(\gamma)$ yields
\begin{equation}
\cosh(\gamma) = 1 + 3i\beta + \frac{\alpha^2}{2} - \frac{3\beta^2}{2} - \frac{3i\beta^3}{10} + \mathcal O(\beta^4) + \mathcal O(\alpha^{4}). 
    \label{C5}
\end{equation}
Note that the coefficients associated with the $n^{\rm th}$ exponents of $\beta$ scale with $6^n/(2n)!$, and thus rapidly go to zero. 
Truncating $\beta$ series in Eq. \eqref{C5} and plugging the result into Eq. \eqref{omegahe} after expanding $\cosh(\alpha)$ and $\sinh(\alpha)$ produces
\begin{equation}
\hat \omega = \left(1 + \frac{i\beta}{2} - \frac{\beta^2}{10}\right)\omega  + \mathcal O(\alpha^{2}), \;\;\; {\rm for} \;\;\; \alpha\ll 1. 
    \label{C6}
\end{equation}

To merge the two expansions in Eqs. \eqref{C3} and \eqref{C6} and approximate $\alpha$ at intermediate values, we can start from Eq. \eqref{C6} and modify the denominator of the second and third term by adding a factor that scales with $\alpha$ and $\alpha^2$, respectively. 
Our numerical experiment shows that the following formula, which matches the above asymptotic expansions, provides a good fit to the original expression in Eq. \eqref{omegah}. 
\begin{equation}
\hat \omega = \left(1 + \frac{3i\beta}{2\alpha \sqrt{1+9\alpha^{-2}}} + \frac{3(i\beta)^2}{2\alpha^2(1+15\alpha^{-2})}\right)\omega. 
    \label{C7}
\end{equation}
Interestingly, the denominators of the second and third term in Eq. \eqref{C7} can be expressed exactly and approximately, respectively, in terms of $\tau$ from Eq. \eqref{tau_supg} (as argued below, the third term has a much smaller contribution than the first and second term, and hence this approximation has little effect on the results). 
Doing so yields
\begin{equation}
\hat \omega = \left(1 + \left(\frac{3i\beta a \tau }{\alpha h}\right) + \frac{1}{2}\left(\frac{3i\beta a \tau}{\alpha h}\right)^2 \right)\omega. 
    \label{C8}
\end{equation}
From Eqs. \eqref{alpha} and \eqref{beta}, $(3\beta a)/(\alpha h) = \omega$. 
Hence Eq. \eqref{C8} can be further simplified to 
\begin{equation}
\hat \omega = \left(1 + i\omega \tau + \frac{1}{2}(i\omega \tau)^2 \right)\omega. 
    \label{C9}
\end{equation}

From the above asymptotic expansion, it is evident that the third term on the right-hand side of Eq. \eqref{C9} is much smaller than the second, as is the second term in comparison to the first term. 
To verify this, one can consider the product of $\tau$ and $\omega$ that scales with $\beta/\alpha$ for large $\alpha$ and with $\beta$ for small values of $\alpha$. 
Since $\beta \lessapprox 1$, in either case $\omega \tau \lessapprox 1$. 
Following this argument, one can include the higher exponents of $\omega \tau$ in Eq. \eqref{C9} with prefactors that quickly go to zero and write as  
\begin{equation}
\hat \omega = \left(\sum_{n=0}^\infty \frac{1}{n!}(i\omega \tau)^n \right)\omega. 
    \label{C10}
\end{equation}
The series in Eq. \eqref{C10} is the Taylor expansion of $\exp(i\omega\tau)$. 
Replacing it will produce Eq. \eqref{omegah}, thus completing this derivation.

To further verify that this approximation is indeed a good approximation of the exact expression for $\hat \omega$, the readers are referred to Figure \ref{fig:tauo}. 

Lastly, we should note that even though Eq. \eqref{omegah} was derived for $\beta \lessapprox 1$, it produces reasonable approximation up to $\beta <5$. 
The values of $\beta$ larger than 5 have little relevance as those signify extremely oscillatory regimes on a very coarse grid that can not be resolved even if the exact form of $\hat \omega$ is employed in higher dimensions. 
In fact, in such scenarios, it is favorable to use the approximate $\hat \omega$ rather than its exact form as Eq. \eqref{omegahe} exhibits an erratic behavior that creates a numerical stability issue. 
That is in contrast to Eq. \eqref{omegah} which behaves well, even at relatively large values of $\beta$.

\section{Consistency analysis} 
\label{app:cons}
Our goal here is to demonstrate the discrete forms discussed in Eqs. \eqref{1D-weak}, \eqref{1D-SUPG}, \eqref{1D-vms}, and \eqref{1D-ASU} will recover the original differential equation in Eq. \eqref{1D-fourier} as the mesh size $h\to 0$. 
Our general strategy is to analyze the nodal solution found in Eq. \eqref{A9} for the Galerkin's method to assess how error changes as $h\to0$. 
To generalize our findings to other stabilized methods, we rely on the modified forms of parameters in Table \ref{table:sum}.

Recall from Eq. \eqref{A9} that
\begin{equation}
    \left(\frac{i\omega h}{6} + \frac{a}{2} -\frac{k}{h}\right)\phi(x_{A+1}) + \left(\frac{2i\omega h}{3} + \frac{2\kappa}{h}\right)\phi(x_A) + \left( \frac{i\omega h}{6} - \frac{a}{2} -\frac{k}{h} \right)\phi(x_{A-1}) = 0,
    \label{E1}
\end{equation}
where $\phi(x_A) = U_A$ is the solution at node $A$. 
For uniform elements of size $h$, the solution at nodes $A\pm 1$ can be expressed in terms of solution and its derivatives at node $A$ using the Taylor series expansion. The result is
\begin{equation}
    \phi(x_{A\pm1}) = \phi(x_A) \pm \phi_{,x}(x_A) h + \phi_{,xx}(x_A) \frac{h^2}{2} \pm \phi_{,xxx}(x_A) \frac{h^3}{6} + \phi_{,xxxx}(x_A) \frac{h^4}{24} + O(h^5). 
    \label{E2}
\end{equation}
Combining Eqs. \eqref{E1} and \eqref{E2} and simplifying the result produces the following equation at $x=x_A$
\begin{equation}
    i\omega \phi + a  \phi_{,x} - \kappa  \phi_{,xx} + \frac{h^2}{12}\left(i\omega \phi + 2a \phi_{,x} - \kappa \phi_{,xx} \right)_{,xx} + O(h^4)=0. 
    \label{E3}
\end{equation}
From this relationship, it is clear that one will recover the original differential equation in Eq. \eqref{1D-fourier} as $h\to 0$, proving the consistency of the Galerkin's method and its second-order accuracy in the diffusive limit. 
One could arrive at the same result by showing the equivalency of the weak and strong form and also the fact the solution to the discrete form approaches that of the weak form as $h\to 0$.

Provided the result in Eq. \eqref{E3}, it is simple to show the consistency of the remaining stabilized methods. Note that Eq. \eqref{E3} also holds for other stabilized methods if one replaces $\omega$, $a$, and $\kappa$ by $\hat \omega$, $\hat a$, and $\hat \kappa$, respectively, according to Table \ref{table:sum} for a given method. Taking the SUPG method, for example, we have 
\begin{equation}
    i\omega \phi + (1-i\omega\tau)a  \phi_{,x} - (\kappa + a^2\tau)  \phi_{,xx} + O(h^2)=0. 
    \label{E4}
\end{equation}
Additionally, note that from Eqs. \eqref{tau_conv} and \eqref{tau_diff} that $\tau^{-2}_{\rm diff} \gg \tau^{-2}_{\rm conv}$ as $h\to 0$. Thus, from Eq. \eqref{tau_supg}, in this limit $\tau \approx \tau_{\rm diff} \propto h^2$. Therefore, Eq. \eqref{E4} simplifies to 
\begin{equation}
    i\omega \phi + a  \phi_{,x} - \kappa  \phi_{,xx} + O(h^2)=0, 
    \label{E5}
\end{equation}
which indicates the SUPG is also consistent and at least a second-order method as $h\to 0$.

Using the same logic, it is straightforward to show that the remaining methods are also consistent and at least second-order accurate.

\section{Derivation of Eq. \eqref{2d-sol}} 
\label{sec:2dsol}
Our overall approach to solving the 2D convection-diffusion problem in the square-shaped domain (Eq. \eqref{2d-prob}) is to first homogenize the boundary condition on the top, then solve the resulting PDE through the method of separation of variables. 
Namely, we are seeking a solution with the form
\begin{equation}
    \phi(x,y) = \psi(x,y) + V(y),
    \label{D1}
\end{equation}
so that $\psi(x,y)$ is governed by 
\begin{equation}
\begin{split}
i\omega \psi + a\psi_{,x} & = \kappa (\psi_{,xx} + \psi_{,yy}),\\ 
\psi(x,0) & = \psi(x,L) = 0, \\
\psi(L,y) &= -V(y), \\
\psi(0,y) & = 1- V(y). \\
\end{split}
\label{D2}
\end{equation}

It follows from Eqs. \eqref{2d-prob}, \eqref{D1}, and \eqref{D2} that $V(y)$ must be satisfied by 
\begin{equation}
\begin{split}  
    i\omega V &=\kappa V_{,yy}, \\
    V(0) &= 0, \\
    V(L) &= 1. 
\end{split}
\label{D3}
\end{equation}

Equation \eqref{D3} is a constant coefficient second-order ODE and its solution can be obtained through elementary means.
It can be expressed as 
\begin{equation}
    V(y) = C_1 \sinh\left(\sqrt{\frac{\omega}{\kappa} i} y\right) + C_2 \cosh\left(\sqrt{\frac{\omega}{\kappa} i} y\right).
    \label{D4}
\end{equation}
Applying the boundary conditions and noting $W^2 = \omega L^2/\kappa$ produces
\begin{equation}
    V(y) = \frac{\sinh\left(\sqrt i W \frac{y}{L}\right)}{\sinh(\sqrt i W)}.
    \label{D5}
\end{equation}

Having $V(y)$, we can now solve Eq. \eqref{D2} through the separation of variables. 
That is to assume $\psi(x,y)$ is self-similar and can be expressed as 
\begin{equation}
    \psi(x,y) = X(x)Y(y). 
    \label{D6}
\end{equation}

For Eq. \eqref{D2} to have a solution of the form shown in Eq. \eqref{D6}, we must have 
\begin{equation}
    \frac{i\omega}{\kappa} + \frac{a}{\kappa} \frac{X_{,x}}{X} - \frac{X_{,xx}}{X} = \frac{Y_{,yy}}{Y} = -\lambda^2,
    \label{D7}
\end{equation}
where $\lambda^2$ are the eigenvalues of the $y$-ODE that is 
\begin{equation}
\begin{split}  
    Y_{_,yy} + \lambda^2 Y &= 0, \\
    Y(0) = Y(L) &= 0. 
\end{split}
\label{D8}
\end{equation}

Note that zero and positive eigenvalues were not considered in Eq. \eqref{D7} as they generate a trivial solution for $Y(y)$. 
The solution to the eigenvalue problem in Eq. \eqref{D8} is  
\begin{align}
        \lambda_n &= \frac{n\pi}{L}, \; n = 1,2,3,..., \label{D9} \\
        Y_n(y) &= \sin{\left(\frac{n\pi y}{L}\right)}, \label{D10}
\end{align}
where we dropped the constants in front of the eigen solutions as they can be merged into those of the $X$-ODE later on. 

Having the eigenvalues $\lambda$, we can solve the $X$-ODE from Eq. \eqref{D7} that is 
\begin{equation}
 \kappa X_{,xx}-aX_{,x}-\left(\kappa \lambda_n^2+i\omega\right)X = 0.
 \label{D11}
\end{equation}

Equation \eqref{D11} is another constant coefficient second-order ODE with a solution that can be expressed as 
\begin{equation}
    X_n(x) = A_n\exp{\left(r_{-n}\frac{x}{L}\right)} + B_n\exp{\left(r_{+n}\frac{x}{L}\right)},
    \label{D12}
\end{equation}
where $r_{\pm n}$ are provided in Eq. \eqref{2d-sol}. 

To put all the pieces together, we must combine Eqs. \eqref{D1}, \eqref{D5}, \eqref{D10}, and \eqref{D12} to obtain
\begin{equation}
    \phi(x,y) = \frac{\sinh\left(\sqrt i W \frac{y}{L}\right)}{\sinh(\sqrt i W)} + \sum_{n=1}^\infty \left[ A_n \exp\left(r_{-n}\frac{x}{L}\right) + B_n \exp\left(r_{+n}\frac{x}{L}\right) \right]\sin\left(\frac{n\pi y}{L}\right),
    \label{D13}
\end{equation}
which is the same as the expression provided in Eq. \eqref{2d-sol}. 

The last step of the process is to solve for $A_n$ and $B_n$ so that the remaining two boundary conditions are satisfied.
More specifically, we must have
\begin{equation}
    \phi(L,y) = \frac{\sinh\left(\sqrt i W \frac{y}{L}\right)}{\sinh(\sqrt i W)} + \sum_{n=1}^\infty \left[ A_n \exp\left(r_{-n}\right) + B_n \exp\left(r_{+n}\right) \right]\sin\left(\frac{n\pi y}{L}\right) = 0.
    \label{D14}
\end{equation}
Defining constants
\begin{equation}
    b_n =  A_n \exp\left(r_{-n}\right) + B_n \exp\left(r_{+n}\right),
    \label{D15}
\end{equation}
permits us to express Eq. \eqref{D14} as
\begin{equation}
    \sum_{n=1}^\infty b_n\sin\left(\frac{n\pi y}{L}\right) = -\frac{\sinh\left(\sqrt i W \frac{y}{L}\right)}{\sinh(\sqrt i W)},
    \label{D16}
\end{equation}
in which $b_n$ can be computed through the Fourier transform and is 
\begin{equation}
\begin{split}
        b_n &= \frac{2}{L}\int_0^L \left[ -\sin\left(\frac{n\pi y}{L}\right)\frac{\sinh\left(\sqrt i W \frac{y}{L}\right)}{\sinh(\sqrt i W)}\right] {\rm d}y,\\
            &= \frac{2n\pi \cos(n\pi)}{iW^2+(n\pi)^2}.
\end{split}
\label{D17}
\end{equation}

Similarly, from the other boundary condition 
\begin{equation}
    \phi(0,y) = \frac{\sinh\left(\sqrt i W \frac{y}{L}\right)}{\sinh(\sqrt i W)} + \sum_{n=1}^\infty \left( A_n + B_n \right)\sin\left(\frac{n\pi y}{L}\right) = 1.
    \label{D18}
\end{equation}
This time, defining 
\begin{equation}
    c_n =  A_n + B_n,
    \label{D19}
\end{equation}
produces
\begin{equation}
\begin{split}
        c_n &= \frac{2}{L}\int_0^L \left[ \sin\left(\frac{n\pi y}{L}\right)  -\sin\left(\frac{n\pi y}{L}\right) \frac{\sinh\left(\sqrt i W \frac{y}{L}\right)} {\sinh(\sqrt i W)}\right] {\rm d}y, \\
            &= a_n +b_n,
\end{split}
\label{D20}
\end{equation}
where
\begin{equation}
    a_n = \frac{2(1-\cos(n\pi))}{n\pi}.
    \label{D21}
\end{equation}

$A_n$ and $B_n$ are obtained using Eqs. \eqref{D15}, \eqref{D19}, and \eqref{D20}. The result is what is provided in Eq. \eqref{2d-sol}. That completes this proof. 

\end{document}